# Fractal Geography of the Riemann Zeta Function

Chris King 28 March 2011 v1.7
http://arxiv.org/abs/1103.5274
Web site: http://dhushara.com/geozeta/
Emeritus, Mathematics Department, University of Auckland

**Abstract:** *The quadratic Mandelbrot set has been referred to as the most complex and beautiful object in mathematics and the Riemann Zeta function takes the prize for the most complicated and enigmatic function. Here we elucidate the spectrum of Mandelbrot and Julia sets of Zeta, to unearth the geography of its chaotic and fractal diversities, combining these two extremes into one intrepid journey into the deepest abyss of complex function space.*

**Introduction:**
This paper completes a discovery process I began in 2009, using computational applications I had developed, looking at the 'dark hearts' [1] - the Mandelbrot parameter planes - of a wide variety of complex functions, including the zeta function, to explore the world of complex functions as widely as possible and elucidate universal properties. This year, as I began to re-explore the parameter planes, using a more versatile second generation version of the application, I became literally sucked into the zeta abyss by an unending stream of intriguing new and surprising features, which rapidly grew to the point where I realized I was dealing with an entire geography of complex function space, spread before me, vast and diverse, like the continents of Europe and Asia combined. These are, as far as I know, hitherto unexplored, apart from Woon's 1998 paper [2] setting out a basic description of the Julia set of zero and the outlines of the Mandelbrot set as in fig 1.

The current paper provides a full investigation of the dynamics emerging from all types of critical point, from those on the real line to the ones adjacent to the critical line x = ½. The software to perform these investigations consists of an open source XCode application for Mac downloadable from: http://dhushara.com/DarkHeart/. For those unfamiliar with complex numbers, discrete dynamics, or the zeta function, there is an extended introduction at the end of the paper.

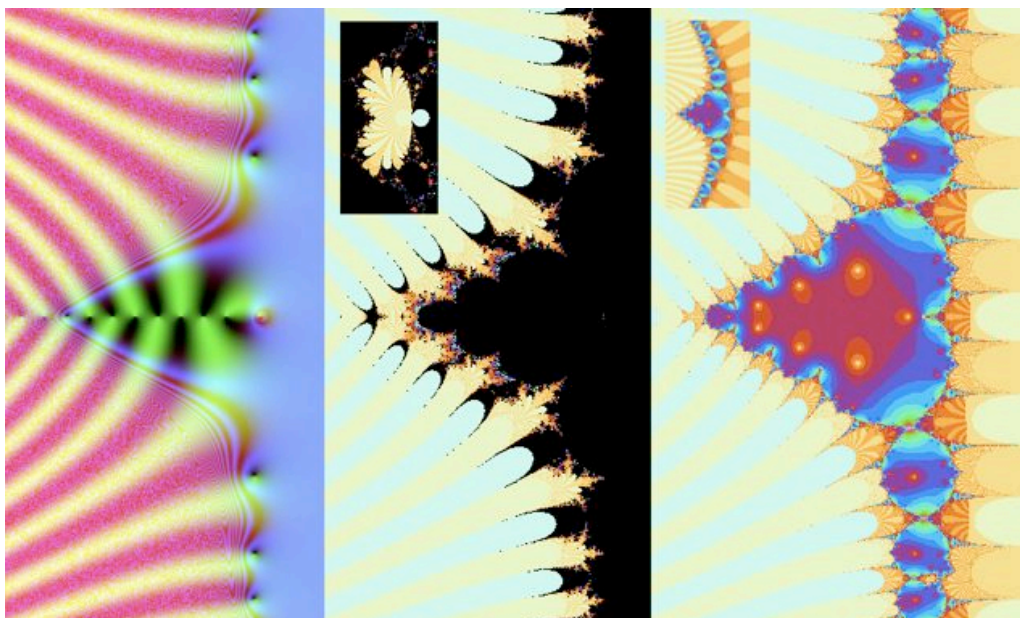

Fig 1: (Left) The zeta function $\zeta(z)$ parametrized by additive colors angle (green/yellow) and amplitude (blue waves and red) so that 0 is black/green, 1 is blue and large values are red and yellow with waves of blue. (Centre) Parameter plane of $\zeta(z) + c$ from the quasi-critical point 1000 on the asymptotic plateau in the right half-plane with singularity 'island' (inset) connected by a fractal thread. The bands of blue and yellow distinguish points iterating towards $\pm\infty$. Strictly speaking the green areas should also be black, as they iterate to the positive half-plane and become fixed, although far outside the iteration limit of the method. (Right) the Julia set of $\zeta(z) + 0$ bounds basins of attraction to the fixed point $\alpha =\sim -0.2959$, containing the non-trivial zeros. Smaller island replicates of the main connected component surround successive trivial zeros (inset). The frond spacing as we ascend to larger imaginary values is spaced irregularly with the zeta zeroes.

We are going to use graphic imagery to explore and confirm the dynamics. The approach is unashamedly numerical, depending on finite approximations of arbitrary accuracy, using computational algorithms. It is also intentionally Zen in its mathematical approach - 'symbolically silent' in its primary use of graphical representations, with minimal symbolic abstraction, to elucidate the geography as fully as possible before describing it. This, supported by the software design is qualitative mathematics in action, using the 'art' to establish the 'math'.

Since we know from Galois that we can't solve fifth degree polynomials symbolically, let alone equations involving transcendental functions, and the Riemann hypothesis that zeta's non-real zeros all lie on the critical line $x = ½$ remains unsolved, despite having been proved for more abstract systems [3], even though all such zeros of zeta are palpably on the critical line, it is clear that the symbolic approach, despite its capacity for abstract generalization, has limitations when dealing with irregular systems of infinite complexity. Hence the research approach taken in this paper.

**A Bridge over Turbulent Waters**

The Zeta function is defined for $real(z) > 1$, either as a sum over powers of the integers, or as a product over primes $\zeta(z) = \sum_{n=1}^{\infty} n^{-z} = \prod_{p\ prime} (1 - p^{-z})^{-1}$. The sum formula is extended to $real(z) > 0$ by expressing it in terms of the eta function's alternating series $\zeta(z) = (1 - 2^{1-z})^{-1} \sum_{n=1}^{\infty} (-1)^{n+1} n^{-z}$. It is then extended again by analytic continuation to $real(z) \leq 0$ $\zeta(z) = 2^z \pi^{-1+z} \sin\left(\frac{\pi z}{2}\right) \Gamma(1-z) \zeta(1-z)$, where $\Gamma(z) = \int_0^{\infty} t^{z-1} e^{-t} dt$, is the gamma function, generalizing the integer factorial $n!$ The result is the most complicated enigmatic complex function known to the human mind.

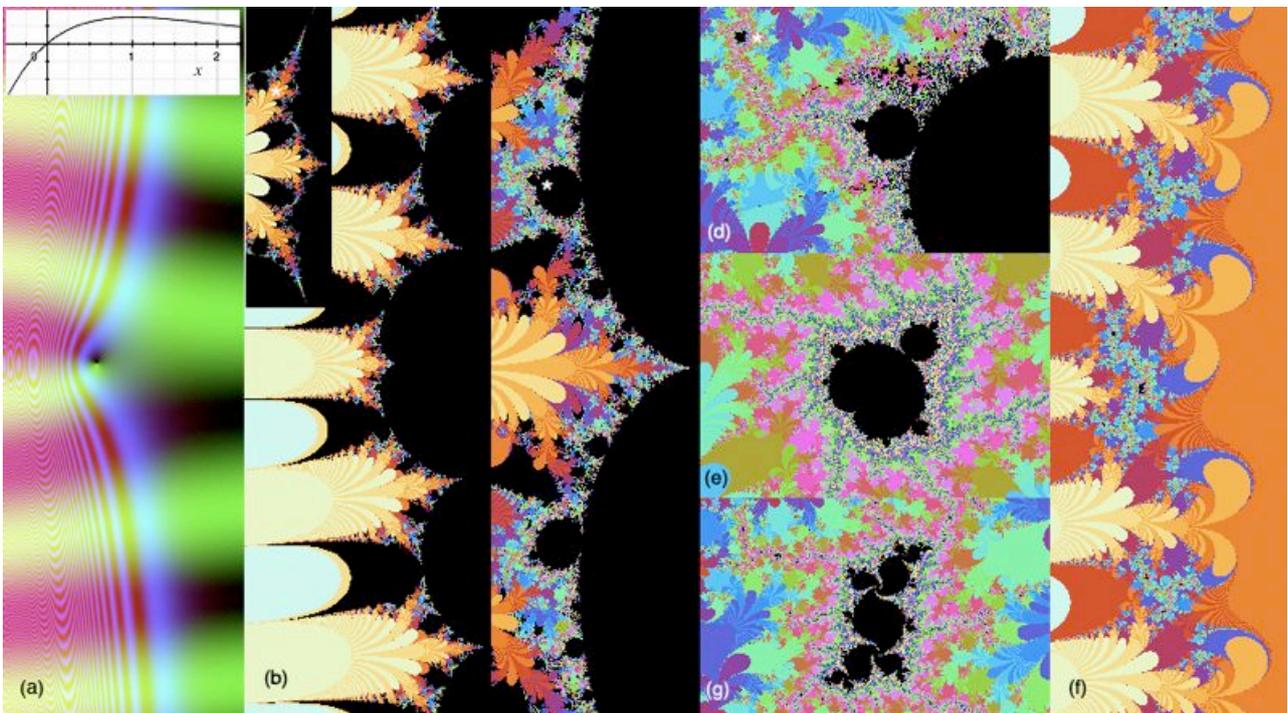

Fig 2: (a) $f(z) = ze^{-z}$ as a real and complex function. (b) Additive Mandelbrot set of $f_c(z) = ze^{-z} + c$ with complex exponential fronds. (c) The central frond for critical point $z=1$ has local quadratic bulbs. The corresponding view for asymptotic plateau quasi-critical $z=1000$ is in inset (b). The bulbs (d) have dendrites supporting quadratic Mandelbrot satellites (e) whose left period 3 bulb has a quadratic period 3 Julia kernel (f,g) and orange plateau matching * inset (b).

To make a transition to the perplexing situation posed by the extreme complexity of the zeta function, let us look at a function that displays pivotal features of the situation in a simpler form. Consider $f(z) = ze^{-z}$. This is an exponential function with an extra $z$ term, which gives it a critical point at $z = 1$, since $f'(z) = (z-1)e^{-z} = 0$ for $z = 1$. All transcendental functions can be represented as power series equivalent to an infinite polynomial. $f(z) = ze^{-z} = z\left(1 + z + \cdots + \dfrac{z^n}{n!}\right) = \sum_{n=1}^{\infty} \dfrac{z^{n+1}}{n!}$.

Every fully differentiable 'analytic' complex function can be so represented. This is similar to zeta but different in an important way. A power series consists of polynomial terms, fixed integer powers of $z$, but a Dirichlet series like zeta consists of a spectrum of exponential functions of integers. The situation is reversed, with weird and wonderful consequences.

However our exponential $f(z) = ze^{-z}$ does form a Rosetta stone for zeta's dynamics. In fig 2 is shown the function, the Mandelbrot set of $f_c(z) = ze^{-z} + c$ from the critical point $z = 1$ (see the end section if this is unfamiliar to you) and a period 3 Julia set. The function tends to zero in the right half plane and to infinity in the left. However, it is sinusoidal in the imaginary direction because an exponential $e^{iy}$ of imaginary $y$ is $\sin(y) + i\cos(y)$, a sinusoidal function whose angle varies with $y$, neatly making complex exponential and trigonometric functions imaginary versions of one another. Notice also the dimple at zero, indicating $f(0) = 0$ - the one zero of the function.

Looking at the additive Mandelbrot set $M_1$ of $f_c(z) = ze^{-z} + c$ from 1, we see it is similar to our unfamiliar zeta case, with exponential fronds representing the waves of the imaginary exponential, zoomed laterally by the real exponential. Now the central frond looks a little different. When we zoom in on it (c), we find it has bulbs just like the quadratic Mandelbrot of fig 35, and these bulbs lead to dendrites containing satellite Mandelbrots identical to the quadratic case, as expressed in Douady and Hubbard's seminal article on polynomial-like mappings [4]. Moreover, when we look at the Julia set of the left-hand 3 bulb on the above satellite, the Julia set (f, g) has a tiny period 3 quadratic Julia 'kernel', set in a fractal web connected to other like kernels.

Now there is another 'quasi-critical' point of this function, where it tends to 0 at (+)infinity. If we had instead used 1000 as our starting point, we would have found a subtly different Mandelbrot set, $M_\infty$ with no quadratic bulbs (inset fig 3(b)), which at the same point as our little Mandelbrot satellite was in the middle of an orange tongue (*). The Mandelbrot set $M_\infty$ classifies the dynamics in the positive half plane, while $M_1$ describes the local polynomial dynamics in the Julia set, as can be seen in fig 2. Julia set dynamics is thus regionally defined in terms of two distinct critical points.

The individual functions in the zeta sum are integer exponentials $f(z) = n^{-z}$ looking like $f$ except for the absence of the forking at zero caused by the $z$ term, having imaginary wavelengths varying logarithmically with $n$, since $n^{-z} = e^{-\ln(n)z} = e^{-\ln(n)(x+iy)} = n^{-x}(\cos(\ln(n)y) + i\sin(\ln(n)y)$. It is the overlapping of these wave functions which gives rise to the irregular pattern of the zeros on $x = \frac{1}{2}$.

**Chasing the Critical Points and their Parameter Planes**

To understand the complex dynamics of zeta we need to examine its critical points. These are precisely the zeros of the derivative of zeta, whose $z$ values are the slope of zeta at $z$, as shown in the right of fig 3. Just as zeta has so-called 'trivial' zeros along $y < 0$ and 'non-trivial' ones on the critical line $x = \frac{1}{2}$, the critical points of zeta are of the same two divergent types, which I will term 'real' and 'unreal', one series along the negative real axis and the other close to the critical line.

Fig 3 shows the first few critical points on $y = 0$ lying between the trivial zeros, and those in the complex plane lying between the non-trivial zeros. The 'real' criticals have oscillating values forming an exponentially varying series $\zeta(-x) \to -2^{-x}\pi^{-1-x}\sin\left(\frac{\pi x}{2}\right)\Gamma(1+x)$. The 'unreal' ones have similar critical values to one another, irregularly wandering between 0.4 and 1. We will name the criticals by rounding down, so the reals we consider are $z$-2, $z$-4, $z$-7, $z$-9, $z$-11, $z$-13, $z$-15, $z$-17 etc. Notice that the 'miniscule criticals' up to $z$-13 lie in the central valley of dzeta where the absolute derivative is less than 1, with $z$-15 forming a transition point and the 'vast' ones, from $z$-17 on, are lost in tiny pockets in the exponentiating highlands. We name the 'unreals' looking at their positive imaginary values e.g. $z$23, $z$31, $z$38, $z$42, … $z$65. They also tend to be located in regions where absolute dzeta is less than 1.

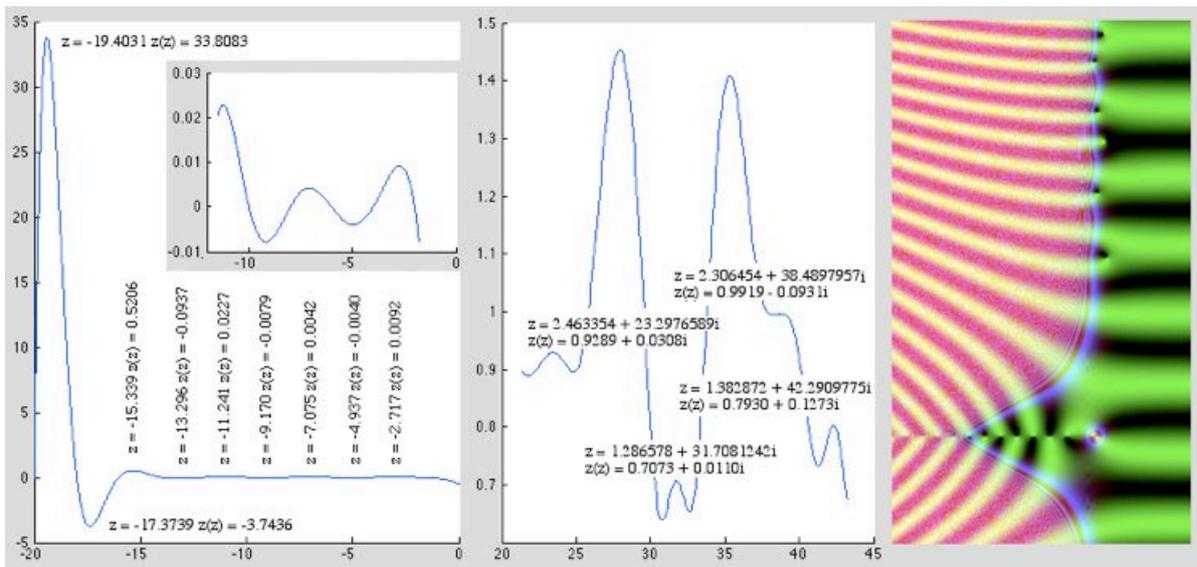

Fig 3: The 'real' critical points of zeta lying along the $x$-axis (left) and the 'unreal' ones close to the critical line (centre). The derivative of zeta (right) shows the two series of critical points as its zeros at the nipples and dimples along the negative real axis and vertically along the outer edge of the blue curve where dzeta has absolute value 1.

However, we are not just looking for critical points, but the places where critical points might iterate to, Mandelbrot satellites that classify interesting Julia set dynamics, which might be somewhere else all together. Looking for tiny regions in a complex exponential fractal can be worse than trying to find a needle in a haystack, so we need at least minimal GPS navigation.

In the quadratic case $f(z) = z^2 + c$ (see end section), the critical point at zero iterates $0 \to 0 + c = c \to c^2 + c$. For $c \sim 0$ we are in the main cardioid, where all points head to the fixed point 0. We can solve for this fixed point. The simplest case is the critical point itself being fixed $c = c^2 + c$, so $c = 0$. Here the $c$ value turns out to be the same as the critical point, but in general, this $c$ value, which we call the 'principal point', could be different from that of the critical point.

More generally, we can try to solve for $c$ values that become eventually fixed or eventually periodic with period $n$ in $m$ steps. These points are the repelling Misiurewicz points, forming the tips and $n$-connection points of the period $n$ dendrites, as well as the attracting main and satellite Mandelbrot sets, in the quadratic case. We will call these collectively 'M-points'. Solving for fixed critical values, fixed at the second step, gives $(c^2 + c)^2 + c = c^2 + c$ giving $c^4 - 2c^2 = 0$ or $c = 0, -2$. These are our original point and the tip of the dendrite on the negative real axis. If the absolute derivative is less than 1 the point is attracting. We thus need to check what these do, by checking whether

$|f'(z)| = |2z| < 1$. The first is (super)-attracting since its derivative is 0. The second however is repelling, since its derivative is -4. Hence it doesn't lead to a Mandelbrot, but the tip of a dendrite.

For our purposes, we seek the simplest of these solutions for the most horrendous function. We won't be able to solve all the equations but we might be able to get a numerical solution and even one that we can display graphically in a useful form.

The very simplest - the critical point being fixed is $c_p : c_r = \zeta(c_r) + c_p$, $c_p = c_r - \zeta(c_r)$, or in the multiplicative case $c_p : c_r = c_p\zeta(c_r)$, $c_p = c_r / \zeta(c_r)$. Since $c_r$ is critical, its derivative is zero, so it is super-attracting and must have a critical value in the Mandelbrot set or its satellites. These 'principal points' $c_p$ can be far from the critical point, even in regions of dzeta where the values are exponentiating towards the infinite.

If we go one step further and look for $c$ values for which the critical value is a fixed point, call them 'fixed values', for the additive zeta Mandelbrot $\zeta_c(z) = \zeta(z) + c$ as in fig 1, we seek $\zeta(c_r) + c = \zeta(\zeta(c_r) + c) + c$. Solving we get $\zeta(\zeta(c_r) + c) - \zeta(c_r) = 0$, or $\zeta(v_r + c) - v_r = 0$, where $v_r = \zeta(c_r)$ is the critical value. For the multiplicative case $\zeta_c(z) = c\zeta(z)$, we get $\zeta(cv_r) - v_r = 0$.

In both cases these 'transfer functions' of $c$ are just transformed copies of zeta, translated, or scaled, in the domain and raised, or sunken, in the range. We can identify the principal point among the fixed values, because it has a 'double twist' in its angle, leading to two yellow angle rays.

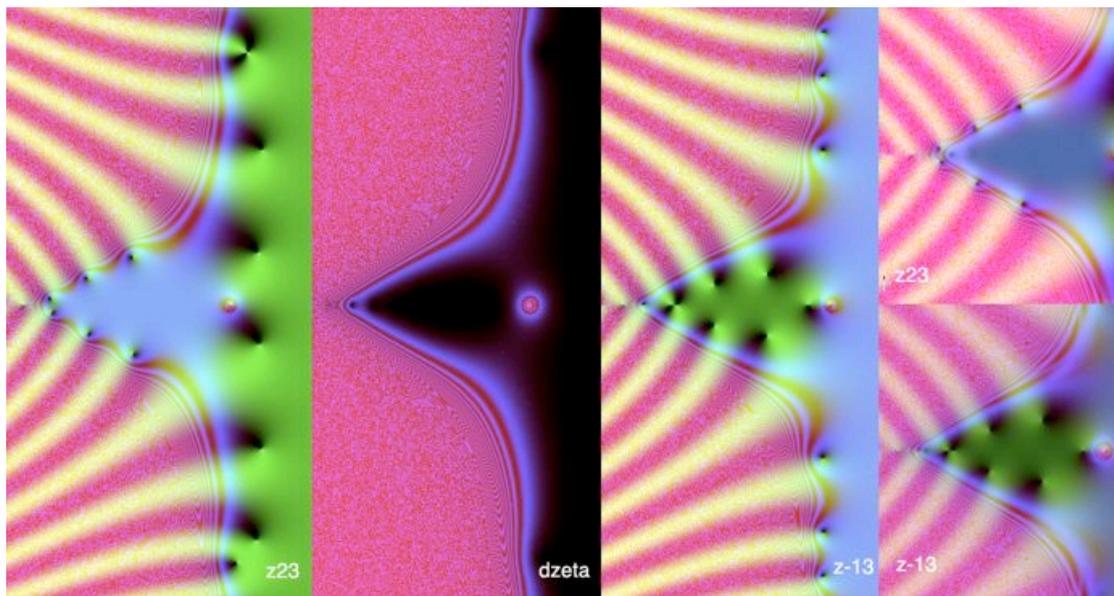

Fig 4: Transformed zeta transfer functions for the critical points $z23$ and $z-13$ compared with that of dzeta with the angle colouring omitted to emphasize the transition across abs($z$) =1. As noted (right) the loci in the bay for $z23$ have absolute derivative greater than 1 so should be Misiurewicz points, while those of $z-13$ are attracting and should lie in the Mandelbrot set, or its satellites. Principal points are identifiable by their 'double rays' top in $z23$ and left in $z-13$.

We can now make a graphical portrait of the transfer functions, locate their zeros and explore the neighbourhood for its local fractal geography. However we also need to know if the fixed points are attracting and thus lie in satellite Mandelbrot sets, or repelling and thus Misiurewicz points, by testing the derivative of zeta. Fortunately we have a serendipitous graphical way to do this, because the colouring scheme for zeta was chosen to highlight absolute value 1, so applied to dzeta, it gives a colour test of the derivative for attracting or repelling. This can be scaled to examine it more closely, or the derivative can be calculated numerically. The method only gives one basic set of candidates, further of which could be found by solving for later fixed or periodic points.

## A: The Additive World

We are first going to explore the geography of the additive Mandelbrot sets $\zeta_c(z) = \zeta(z) + c$ for the various critical points of zeta and how they interact with one another and with the Julia sets they define. Subsequently we will explore the more bizarre dynamics of the multiplicative parameter planes of $\zeta_c(z) = c\zeta(z)$.

## 1: Far East - the Asymptotically-Critical Plateau

We begin with the Mandelbrot set $M_\infty$ in fig 1, originating from the nominal quasi-critical point 1000 on the plateau of zeta converging to 1 in the right half plane. This does not display polynomial bulbs or satellite quadratic Mandelbrot sets but consists of fractal enclosed representations of bounded Mandelbrot regions, interpenetrated by the chaotic escaping set, fractally replicating the exponential ferns, whose global form is illustrated by the anti-Mandelbrot island around the singularity in the inset of fig 1.

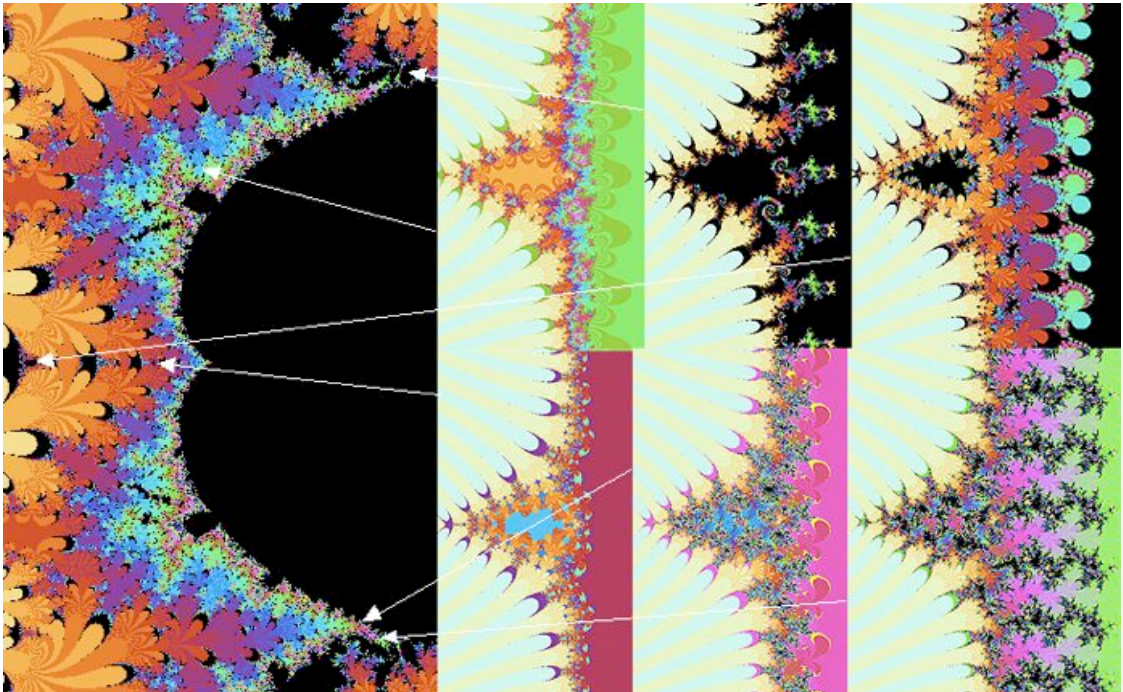

Fig 5: Locations on the Mandelbrot set $M_\infty$ classify the asymptotic dynamics in the right half-plane. The quantitative step-colouring of each location on $M$ coincides closely with the step colouring of dynamical escape on the plateau. Points in $M_\infty$ and its fractal islands remain bound (black).

What the plateau's parameter plane $M_\infty$ is measuring is the dynamics of differing $c$ values, iterated from 1000 far in the right half plane. This is illustrated in fig 5, showing the way atlas addresses on $M_\infty$ define the asymptotic step dynamics in the right half plane. As we pass through fractal regions of $M_\infty$, this results in a fractal sequence of dynamic 'explosions' of the right half-plane whenever a path in $M_\infty$ crosses a boundary between a black and a coloured region, dramatic in movie format.

Points $c$ in the right half-plane iterate to the fixed point $c + 1$ because $\zeta(c) = 1$, so iterating from the quasi-critical point, $1000 \to \zeta(1000) + 1000 = 1001 \to \zeta(1001) + 1000 = 1001$. This is true for all $z$ with large positive real part, so the iteration is fixed. Numerically, the principal point is $1000 - \zeta(1000) = 999$, which again places it in the asymptotic limit, which coincides with the picture of the Mandelbrot set engulfing the positive half-plane.

Although the dynamics consists of fractal exponential fronds, these do display the same mediant-based fractional winding adding fractional rotation periods that the quadratic Mandelbrot bulbs have. In fig 6 we show that the same mediant winding sequences, we see in fig 35 for the quadratic Mandelbrot appear on the bays in the fronds bounding $M$.

Intriguingly the Farey tree of mediants appears in one variant of the Riemann hypothesis. Farey sequences consist of all fractions with denominators up to $n$ in order of magnitude – viz $\frac{0}{1}, \frac{1}{5}, \frac{1}{4}, \frac{1}{3}, \frac{2}{5}, \frac{1}{2}, \frac{3}{5}, \frac{2}{3}, \frac{3}{4}, \frac{4}{5}, \frac{1}{1}$. Each fraction is the mediant of its neighbours $\frac{n_1}{d_1}, \frac{n_2}{d_2} \to \frac{n_1+n_2}{d_1+d_2}$. Two versions of RH state [5]:

$$\text{(i) } \sum_{k=1}^{m_n} |d_{k,n}| = O(n^r), \text{ any } r > 1/2 \text{ and (ii) } \sum_{k=1}^{m_n} d_{k,n}^2 = O(n^r), \text{ any } r > -1$$

$d_{k,n} = a_{k,n} - \frac{k}{m_n}$, where $m_n$ is the length of the Farey sequence $\{a_{k,n}, k = 1, \cdots, m_n\}$

We can colour $M$ according to how many steps it takes to reach within $\varepsilon$ of a fixed point or periodicity and colour by the number of steps in blue and add redness for the period. This immediately shows up the periodicities of the bays neighbouring the boundary, which can be confirmed to correspond to Julia sets with rotational periodicity the same number, confirming the sequences of periods of the bulbs and the mediant relationship in the fractal progression.

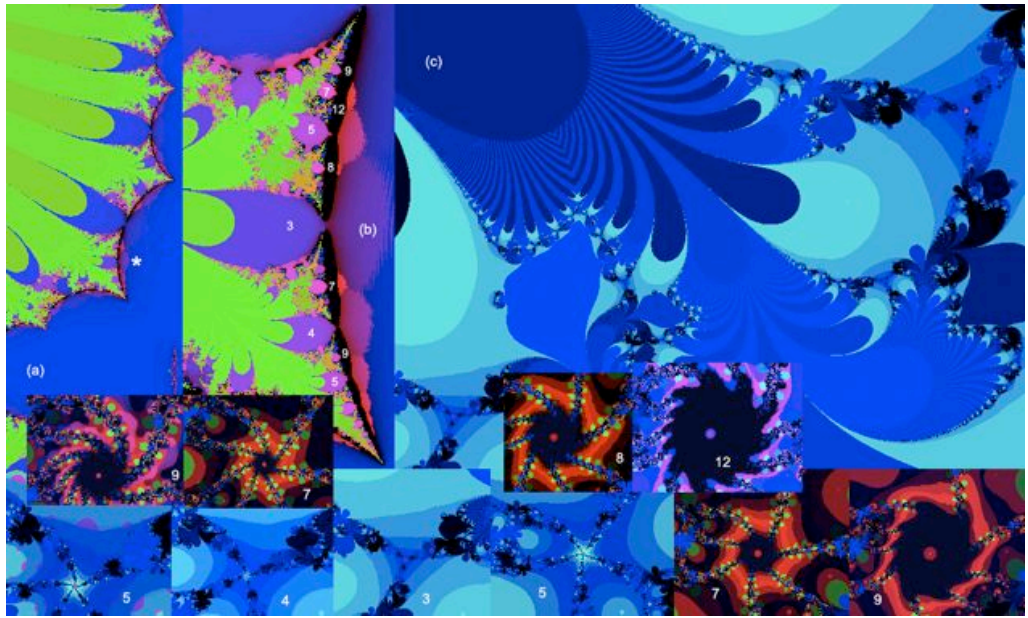

Fig 6: Shading the bulbs demonstrates their periodicity, as confirmed by the Julia set portraits, which display the corresponding rotational periodicities in their spirals, confirming an upward set of odd periodicities 3, 5, 7, 9 ... and a downward set of integer periodicities 3, 4, 5, 6 ..., each with mediant fractality, viz (3,4)=7, (4,5)=9, (3,5)=8, (5,7)=12. Period 3 is shown top right of (c).

## 2: Real Critical Points, from Miniscule to Vast

While the zeta zeros on the real line are regarded as the trivial solutions of $\sin\left(\frac{\pi z}{2}\right) = 0$, the critical points on the $x$-axis are anything but trivial, and each displays qualitative features of zeta that give each critical point a distinct role in the dynamics. When we have a function with more than one critical point, to understand the dynamics, we have to investigate the Mandelbrot set of each critical point. The critical points contribute to different dynamical features of the whole Julia set, as illustrated in figs 5 and 10. The dynamics also involves interactive effects between the critical points which causes their Mandelbrot sets to appear merged or amorphous and the dynamics in different parts of the Julia set to be influenced by each of the critical points.

In this respect the situation is very different from the quadratic case, where the Mandelbrot set is an infinite atlas of the dynamics of the Julia sets, each of which has a single type of dynamics determined by the *c* value of the single critical point. In the case of zeta, with an infinite collection of critical points, the relationship between the Mandelbrot and Julia dynamics is structurally analogous to a Fourier transform. As before, the Mandelbrot set for a given critical value is a spatial 'integral' of Julia dynamics over continuously varying *c* values. However the Julia set dynamics is now determined by a countably-infinite collection of critical points, each of which can fractally dominate the dynamics around its M-points. Examples of Julia set dynamics responding to many critical points are illustrated in several of the figures.

**(a) Continental Divide: The Critical Point *z*-15**

The critical point *z* ~ -15.339 commands a pivotal role in the dynamics of the central basin. When we examine its principal point and fixed values in the central valley, fig 7(b), we find they are in the main body of its Mandelbrot set , close to the shores of the three bays we can also see in the Mandelbrot set of fig 1, originating from the three bounding fronds. Pivotally its principal point is right off the shore of the innermost bay, and it is here we find a sequence of quadratic bulbs and the cusp familiar in the quadratic Mandelbrot set (fig 35).

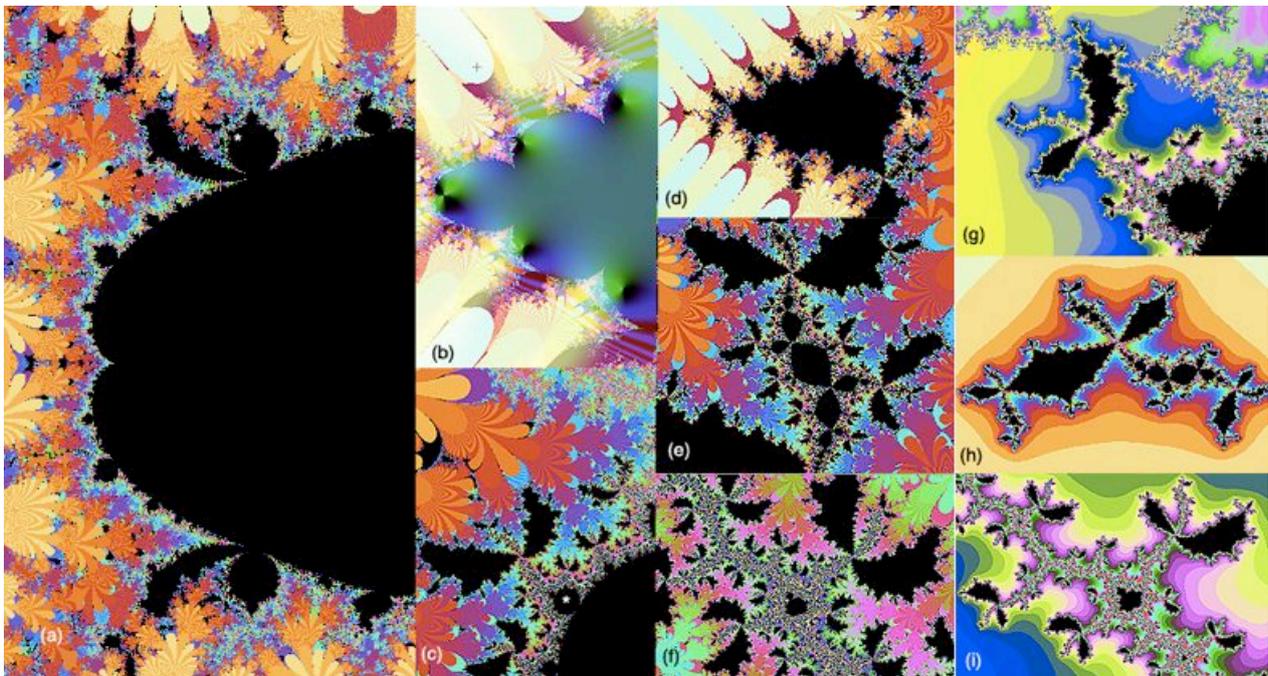

Fig 7: (a) Base of the central valley for *z*-15 showing quadratic bulbs perturbed by cubic and higher dimensional interference. (b) The critical value fixed points, including the principal point, all lie in the central valley close to the boundary, thus dominating the bulb dynamics. (c) The period 3 side bulb gives rise to a Julia set (d,e) with obvious period 3 dynamics. The dynamics are perturbed by adjacent critical points both of which are in a cubic relationship to z-15 and possibly other 'unreal' criticals. The features of (c-e) share dynamical morphology with regions on the Mandelbrot set (g) and the corresponding Julia set (h) of the cubic function $f(z) = z^3 - z + c$. (f, h) Satellite Mandelbrot sets of the two functions also share cubic morphology.

As we move into the cusp, fig 8 lower, we find high periodicity dynamics characteristic of classic quadratic regions such as 'seahorse valley' (6 in fig 35). As we move further away from the cusp the dynamics becomes more complicated, with the largest bulb having an appendage from the base of a kind also seen in cubic functions where the critical maxima and minima are close enough that their dynamics interferes. Although zeta has no degenerate critical points which are multiple zeros of the derivative (compare fig a6), in fig 7 comparison is made between regions of the cubic function $f(z) = z^3 - z + c$ and this region, in terms of both its Mandelbrot and Julia dynamics, confirming the similarities in a Julia set from the period 3 sub-bulb (*) and in satellite Mandelbrots from each. This dynamic interference possibly originates from *z*-13 as it shares features here with *z*-15, however many other critical points could be involved. For example, the unreal critical *z*95 has a deformed version of the *z*-15 structure, which also has the same cubic 'wings'.

Many of the 'unreal' criticals also have critical values close to the critical value of *z*-15 of 0.52 and fixed values in similar locations (see figs 21, 22), so that the dynamics surrounding the central bay consists of the superposition of a countable infinity of perturbations – a little like the humming on telegraph wires in the desert consists in principle of summed vibrations along the transmission line.

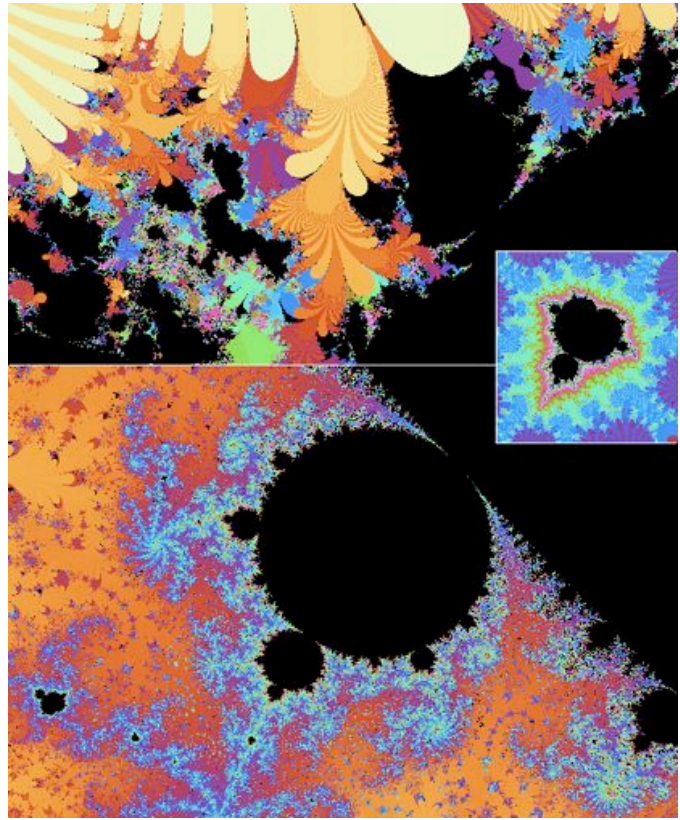

Each frond is a fractal replicate of the entire dynamical parameter plane, so has a fractal replicate of the central valley, increasingly to one side, as we move up successive fronds. For many critical points such as *z*-17, fig 17, the fractal valleys replicate the central valley dynamics, but not in *z*-15, as shown in fig 8 above, where two adjacent generations of fractal valley each have their own distinctive dynamics. Nevertheless the fractal valley in fig 8 does support a Mandelbrot satellite in a corresponding location to the satellite we find in fig 9 at the base of the central valley.

Fig 8: (Above) The structures in fig 7 are not fractally repeated in the fractal valley at y~13 (left) which has distinct dynamics from the central valley, although *z*-15 does have a Mandelbrot island at the starred point (right). Fractal repeats do occur for *z*-17 on (figs 17-19) and the dynamics is more similar for *z*-2 (fig 16). (Below) An exploded view in the cleft of the main basin of *z*-15 showing high periodicity dendrites.

There are also multiple fractal replicates of the central valley interspersed down the base of the valley (see figs 9 and 15) and into the crests and troughs running along the negative real axis, which will be useful in elucidating the dynamics.

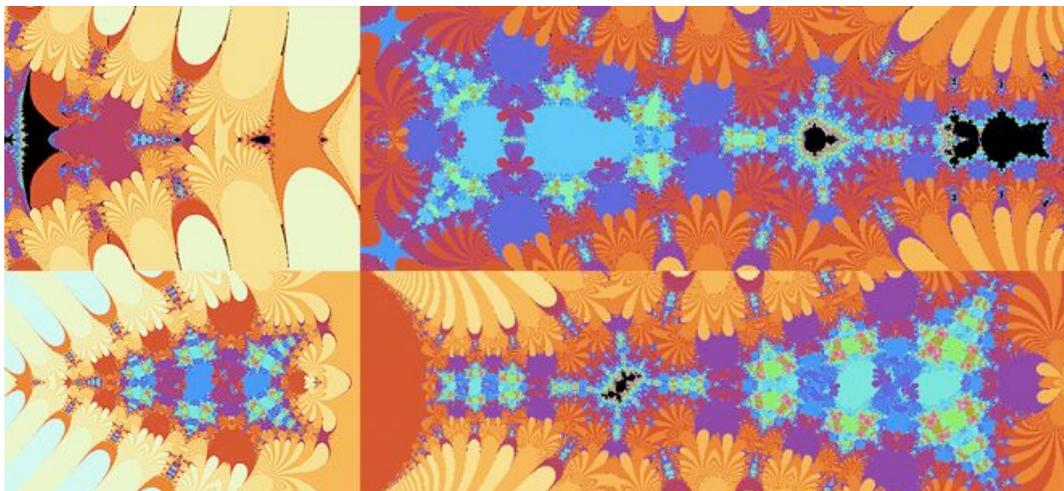

Fig 9: Fractal basin on the real line at *x* ~ -17.95 (top, enlargement right) has a satellite of the critical point at *x* ~ -15. The Julia set of its upper period 3 bulb (below) generates a Julia set with a period 3 kernel in its base.

## (b) Gently Undulating Lowlands : The Miniscule Criticals  $z$-2 – $z$-13

We deal with the miniscule criticals as a group, because, in many ways, they behave like a higher degree polynomial of degree between 4 and 6 depending on the situation.

We start by looking at $z$-2, and $z$-13 at the large bulb we investigated in fig 7. When we examine the miniscules, we find this has become a towering amorphous structure, I will call the 'ant', indicating interference between several critical points. On $z$-13 this has bulbs, indicating these regions are quadratically sensitive to it, which also have satellite Mandelbrot sets on their dendrites (*), as does $z$-15. $z$-2 also has satellites, indicating the 'ant' region is sensitive to most of the miniscules. Notably there is no such structure on $z$1000. In $z$-2 this region is a fractal replicate of the central bay, with three beaches separated by fronds. The 'head' region even has 'horns' which effectively replicate the ant structure in the central bay within itself. When we look at the whole central bay of $z$-2 or $z$-7 we find the 'ant' is a fractal replicate of the bay repeated for each of the three fronds and fractally on all scales and is also present in the bays of the fronds in fig 16. Each of these also has the complex quadratic structure involving several critical points we find in the ant. Period 3 bulbs on each of the satellite Mandelbrots generate period 3 Julia kernels, confirming the satellite of each critical is determining the Julia dynamics in the period 3 web of each set, despite the fact that the Julia set is sensitive only to the location of the $c$ value and not the critical point that generated the satellite. This shows each of the critical points are collectively determining the Julia dynamics.

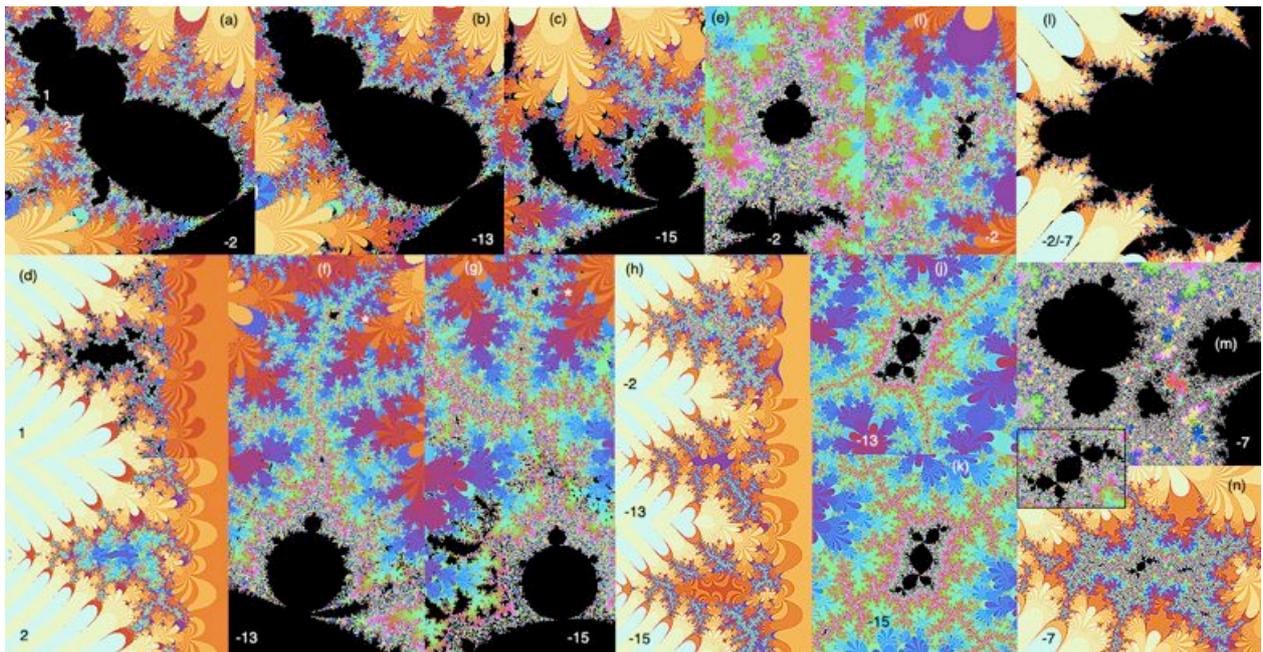

Fig 10:  Critical points from z-2 to z-15 all show fractal polynomial structures on the boundary of the central valley (a–c), including satellite Mandelbrot sets (e-g). These are not only perturbed by the 'miniscule criticals z-2 – z-13 but by many of the unreal critical points, many of which have critical values surrounding that of z-15.  The critical points z-2 to z-9 have critical values very close to 0 and thus form an atlas of the dynamics in the central valley and the zeros.  (d) z-2 classifies the differing central valley dynamics between points 1 and 2 in (a). Period 3 bulb dynamics of the satellite Mandelbrots of the three critical points each show distinct regions of period 3 dynamics in their Julia dynamics (h-k), confirming all three critical points leave their mark on the Julia set. Only that of z-2 continues into the central basin. (l) The central bay for $z$-2/$z$-7 showing the 'ant' is one of three fractal bay structures, which are repeated fractally on all scales. (m) A satellite for $z$-7 at the head of the ant, confirming the head and head cusp are sensitive to the 'miniscules', particularly $z$-2 as shown in fig 11. Its period 3 Julia kernel web (n) is covering the central valley.

The effects of each however differ. $z$-2, and with it, the lesser miniscules, form an atlas of the dynamics passing close 0, as their critical values are very close to 0.  Hence they determine dynamics in the central bay. Sampling the point 1 on the 'ear of the *z*-2 'ant' which lies outside the *z*-13 'ant' gives a connected black centre while the point 2 lying outside all three has a chaotic centre. Notice that the 'ant' is absent in $M_\infty$ and indeed the asymptotic plateau in all the Julia sets is brown indicating escape there. But only in *z*-7 and *z*-2 is the web of the period 3 kernel connected across the central basin. We thus can see in the Julia sets the regional actions of three distinct critical points simultaneously, central basin, asymptotic plateau and local polynomial.

The collective evolution of the miniscules and their undulation in the Mandelbrot sets with their critical values is clearly portrayed in the dynamics of the apex of the innermost bay, where the fixed value in the neighbourhood of *z* = -16 points to the blunt frond apex for all of the miniscules, which undulate in position with their critical values, but becomes a quadratic cusp for *z* = -15, which also dominates the local dynamics of unreal critical points, the asymptotic plateau, and *z*-19. However *z*-2 does have a quadratic cusp with bulbs at the head of the ant and its sibling bays.

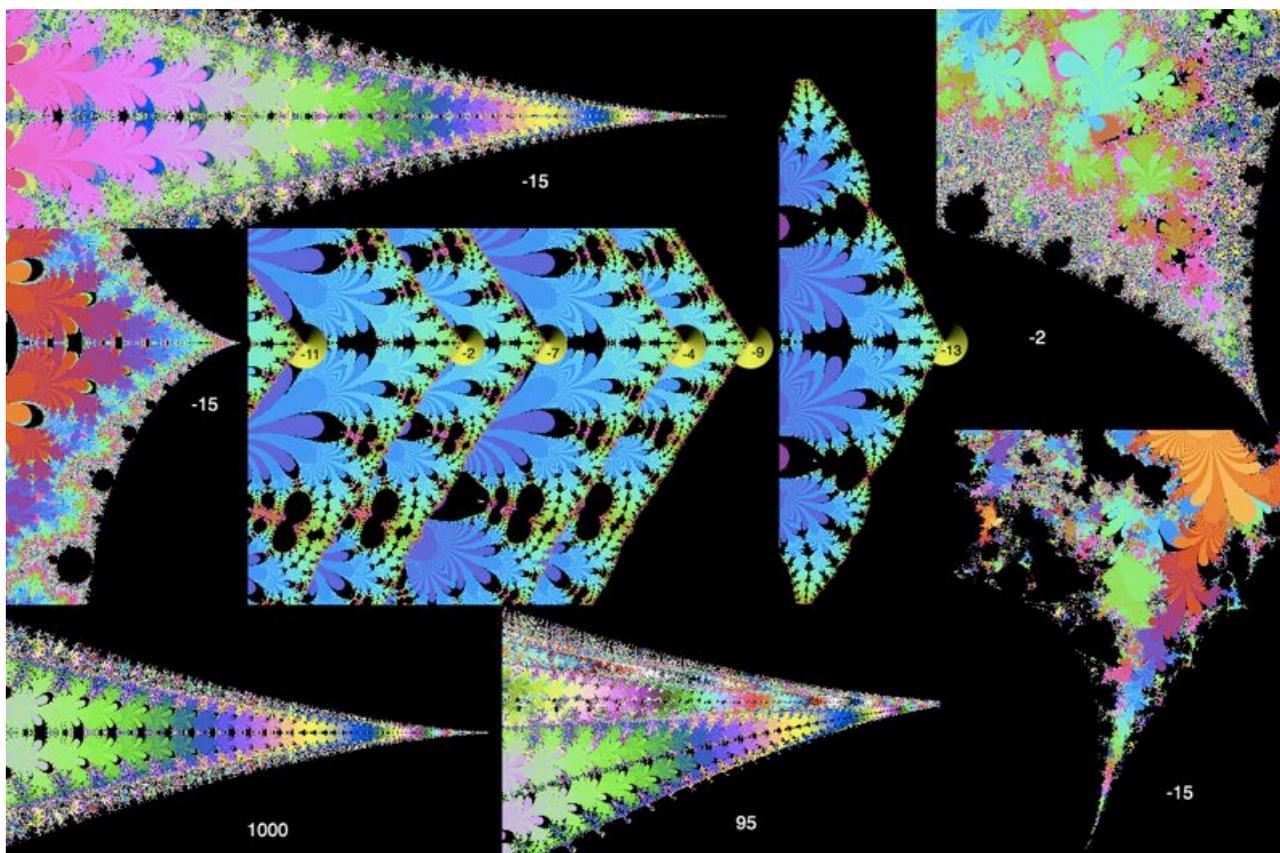

Fig 11: z-15 as critical transition All the critical points from z-2 to z-13 have their fixed value corresponding to the zero at -16 converging to the tip of the basal frond, however when we reach z-15, the entire basin boundary turns into a quadratic cleft. This is conserved by critical points with critical values in the range 0.4-33.8 as illustrated below for *z*1000 and *z*95 and is also true for *z*-19, with a critical value of 33.8, indicating that *z*-15 is influencing the dynamics of all these in this region. At z-17, the central valley becomes flooded (see fig 17). However z-2 has a quadratic cusp on the head of the 'ant' (top right). All these structures differ from the naked cusps at the tip of exponential fronds in the black ocean of the Mandelbrot set (lower right). Images all to scale of 0.02, except for the centre left and right. The transition can be viewed at the website as a video showing interaction between an evolving fractal process interacting with a 'static' representation of the asymptotic limit Mandelbrot set.

We now turn to decoding the collective dynamics of the miniscules and their influence on the dynamics near the real line. Mandelbrot satellites of the miniscules occur in a number of fractal regions n the real line, several of which are fractal replicates of the central valley (see fig 14), which have complex amorphous regions which originate from the overlapping effects of the critical points on one another's dynamics. These regions can be distinguished from a number of fractal regions that are simple basins with a single periodicity, by colouring according to the incipient period. This shows the compound sets have varying periodicities. We can then look for satellite Mandelbrots to confirm they are a Mandelbrot compound structure, as illustrated in fig 12.

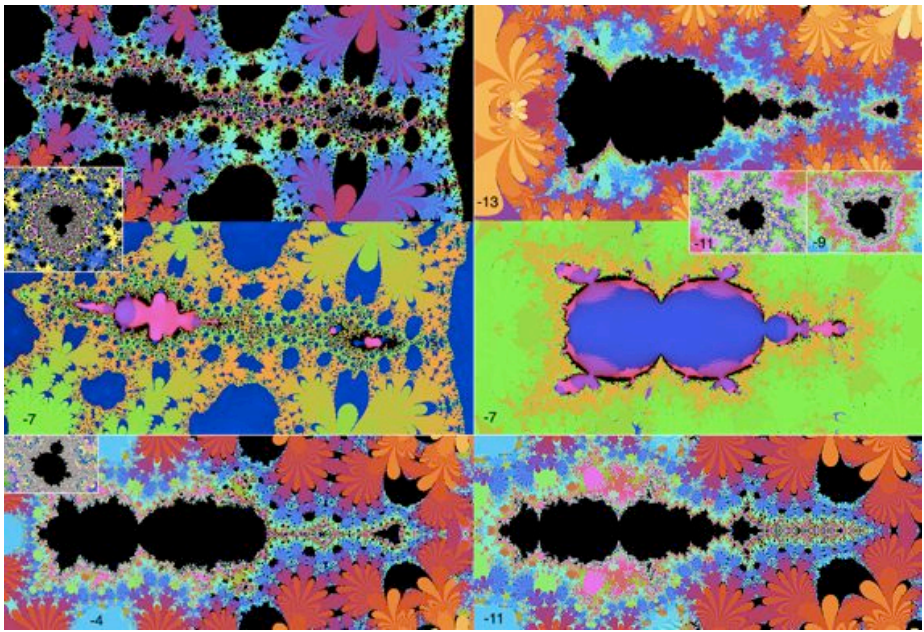

Fig 12: Complex sets displaying overlapping effects of several critical points show their nature through each of them possessing well-formed Mandelbrot satellites, despite having an amorphous morphology. (Top left) region connecting a frond to the central valley. (Top right and below) fractal replicates of the central valley displaying differing degrees of critical point interference (3 and 4 in overview fig 14). Highlighting incipient periodicities (middle row) helps to distinguish complex sets from the blue exponential islands (left) all of which have simple fixed point dynamics.

Figures 12, 13 and 14 show how the relative dynamics of the miniscules can be revealed in stages, by examining each of the regions in positions labeled 1-4 in the top of fig 14. These are each fractal replicates of the central valley and expose the relative dynamics of the miniscules, all of whose principal points are submerged in the central bay.

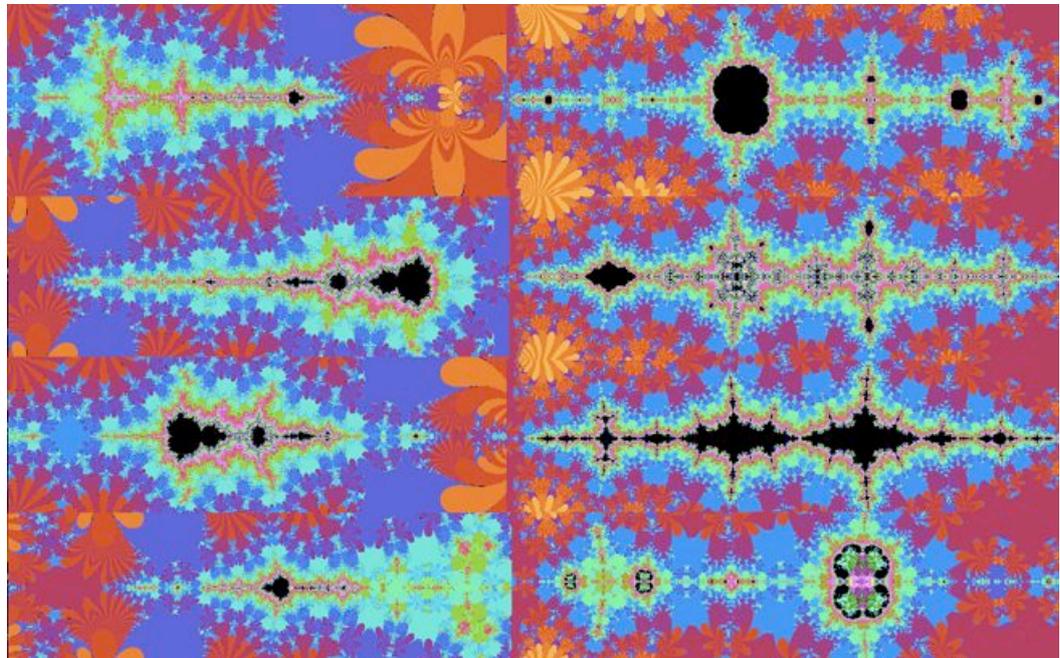

Fig 13: (Left) the local Mandelbrot islands of the first four critical points on the real axis in the fractal replicate 1 in fig 14. (Right) Central valley region of the corresponding Julia sets approximating the $c$ values confirms all four parameter planes influence the Julia dynamics.

The largest at position 3 is the most merged and shows a quadratic satellite only for $z$-13, the most far-flung of the set. The next at position 4 has a greater degree of separation, as shown in fig 12, but still the dynamics is merged, with only $z$-11 showing a clear satellite, despite others having sub-satellites on their surrounding dendrites, confirming this region is a compound Mandelbrot.

Things become much clearer in region 1, where we can see from fig 13 that each of the first four criticals have satellites which are oscillating in position in relation to their critical values, as we saw in fig 11. We lay those of the first four critical points in this region out in sequence, so we can see each as a well-formed 'black heart', each with sensitivities to the location of the 'black hearts' of the other critical points. The evolution of these satellites is confirmed by the varying dynamics of the corresponding Julia sets shown on the right.

When we move to region 2 of fig 14 we find a clear case of fractal separation of the satellites, which now follow a sequence we shall also see extended for all the real criticals in relation to the fronds. Each of the successive criticals forms a graded sequence, with one max-min pair to each of the three frond-pairs, moving from the outermost in the bay to the innermost, laying bare the dynamics of the miniscules, which was submerged in the central bay.

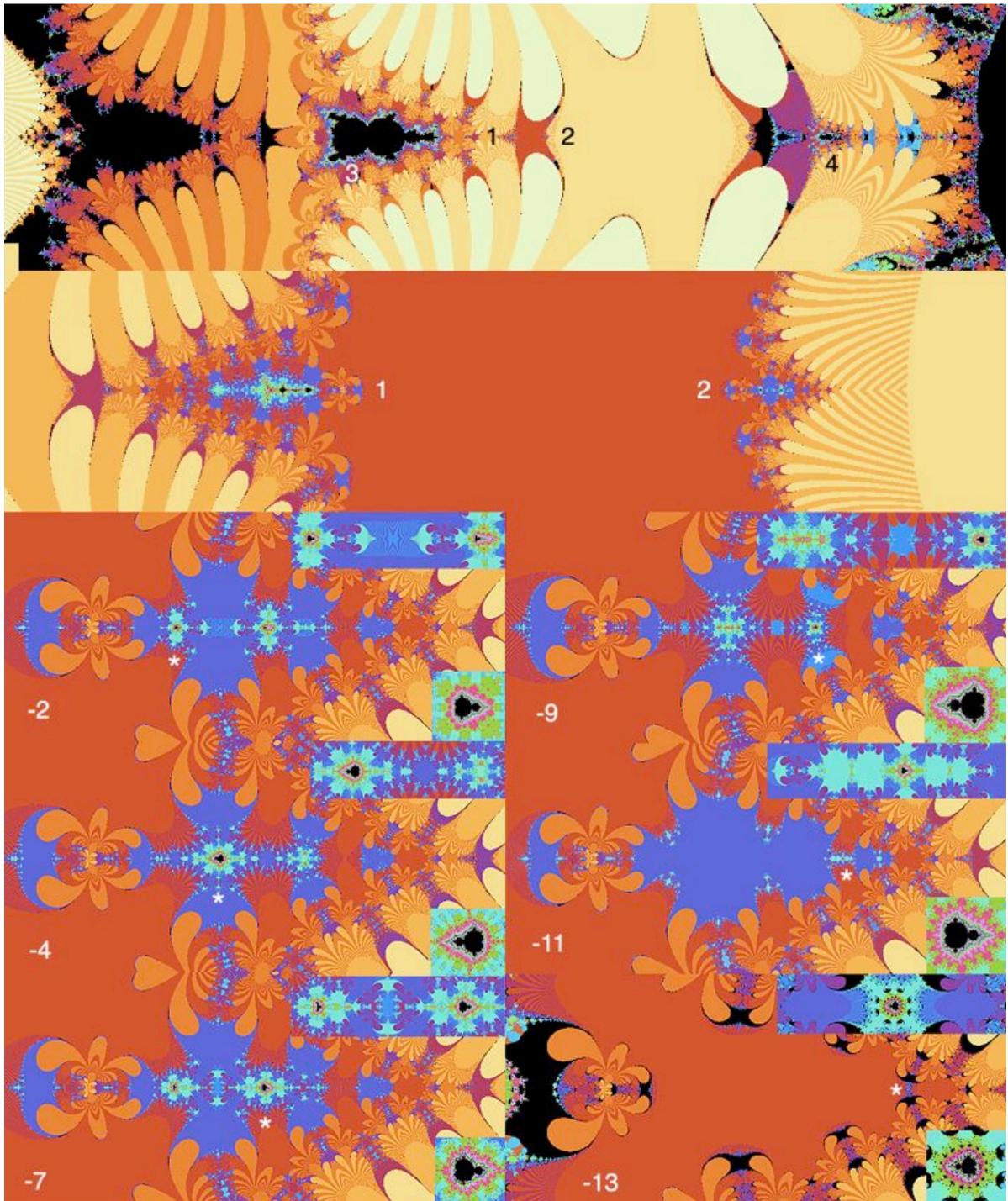

Fig 14: Fractal replicate of the central valley (2) at top shows the evolution of the 'miniscule' critical points, which is otherwise hidden because their fixed values fall into the central bay's blackness. A consistent evolution of the Mandelbrot satellites down the fronds is shown, with one hump and one trough for each successive frond.

The satellites have base periodicity 3 in their central region, when compared against the period 3 satellite on the negative real dendrite of the quadratic of fig 35, as shown in fig 15. Evidence of this can also be seen in the Julia sets, by comparison with a base Julia set for the quadratic satellite. This explains how these satellites can exist in a region where the derivative is large, because for a period 3 cycle we calculate the derivative by the chain rule as a product of the derivatives at the three points in the cycle, one of which is close to zero and has a tiny derivative.

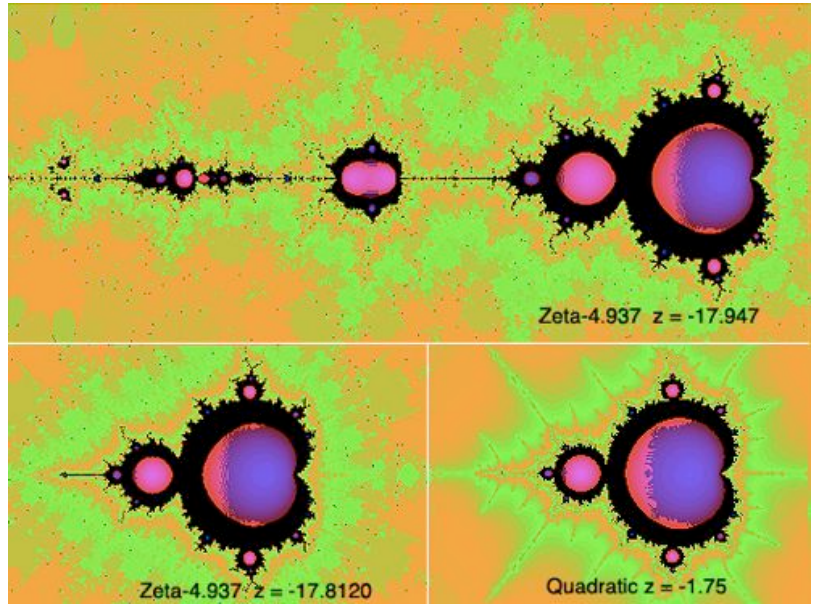

Fig 15: Period-sensitive colouring of Mandelbrot satellites from replicates 1 (top) and 2 (lower left) for z-4 both coincide with the period 3 satellite on the quadratic Mandelbrot set, confirming they are period 3. This both coincides with the forms of the Julia sets in fig 14, which show real period 3 dynamics and explains how they can exist in a region where the derivative has absolute value greater than 1, because other steps in the period 3 cycle include points close to 0 with tiny derivatives, ensuring the period 3 derivative, calculated by multiplying the three derivatives, by the chain rule, confirms the 3-period is attractive overall. For example in the region 2 satellite approximation gives  -17.8120 > -19.8882 > -4.9145 … , with overall derivative -8.8565*101.3019*1.5049e-05 = -0.0135.

This evolution is replicated in the fractal bays present in each frond, as illustrated in fig 16 where there is a homologous evolution in the base of the valley y ~ 20. The dynamics in this valley are very similar to those of y ~ 13, in fig 8, for both *z*-2 and *z*-15.

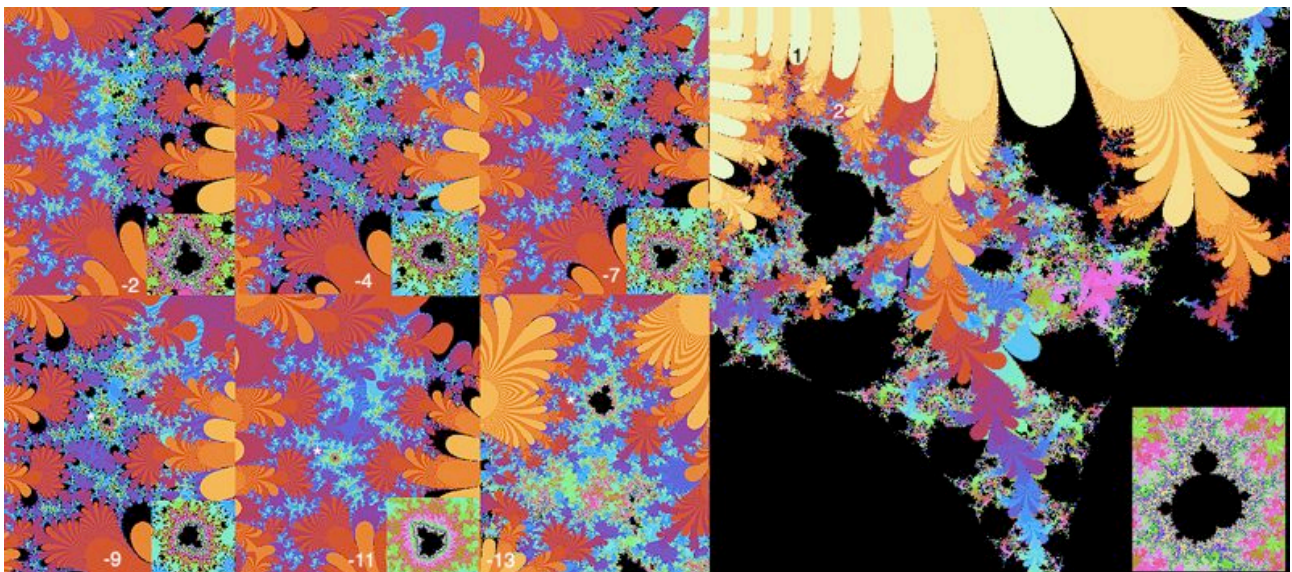

Fig 16: 'A Garden Enclosed is my Beloved' – Base of fractal valley around y ~ 20 (1 right). The evolution of the 'miniscule' critical points z-2 (*) - z-13 (last larger scale) is also presented in the fractal valleys of the fronds. The basin area also supports satellites, as shown in the inset right located at 2 and is a partial homolog of the 'ant'.

**(c) Lofty Peaks of Altiplano – The Vast Criticals**

We now enter a sparse mountainous landscape heading outside the central valley, where the critical values and derivatives become exponentially huge and the transfer function begins to cause large translations, far into the positive and negative reals.

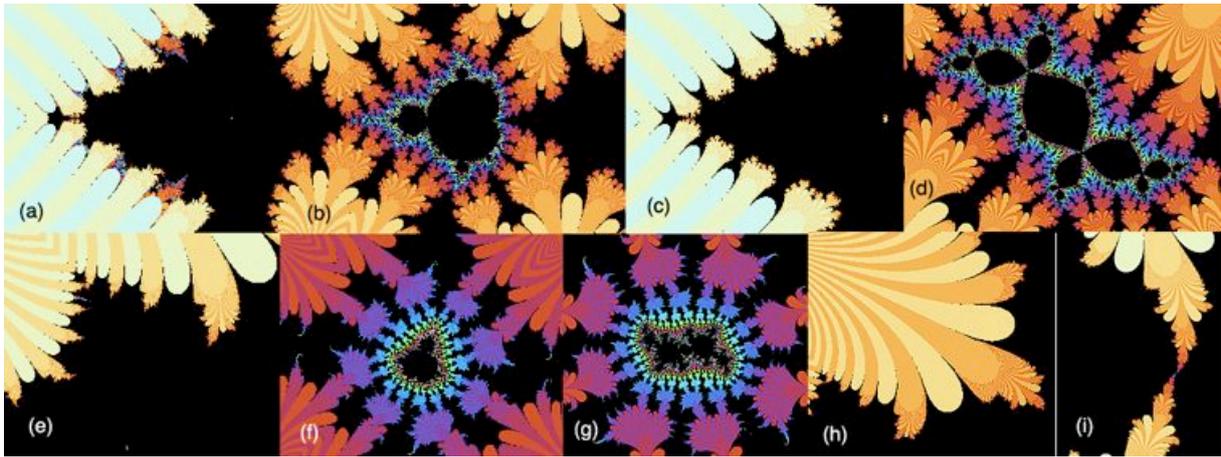

Fig 17: Once we arrive at z-17, the landscape becomes sparse, the central valley becomes submerged and the fronds truncated. The innermost pair of fronds meet in a Mandelbrot set at the principal point (a,b), whose period 3 bulbs generate a Julia set (c,d) with period 3 kernels. The valleys in each frond also have fractal replicates of the Mandelbrot set in the central valley in the same relative position (e,f), Julia set (g), however this is not at the fixed value, which corresponds to the tip of the corresponding frond (h), with Julia set (i) having a touching frond pair. Zeta Misiurewicz points thus include the tips of fronds as well as dendrite *n*-hub points (see figs 20-22), as laso noted in fig 11.

The case of *z*-17 in fig 17 shows the entire central valley flooded back to the fourth frond pair where the two fronds meet in a single Mandelbrot satellite. A *c* value in the period 3 bulb of this gives an equally sparse Julia set with a period 3 Julia kernel held between the same two fronds. This process is fractally replicated in the valley in each frond, with an isolated satellite at the same frond pair. This is however not the location of the fixed values, which lie at the tips of successive fronds and generate Julia sets in which a frond pair are just touching at their tips. This is consistent with the principal point being the only fixed value guaranteed to be attracting.

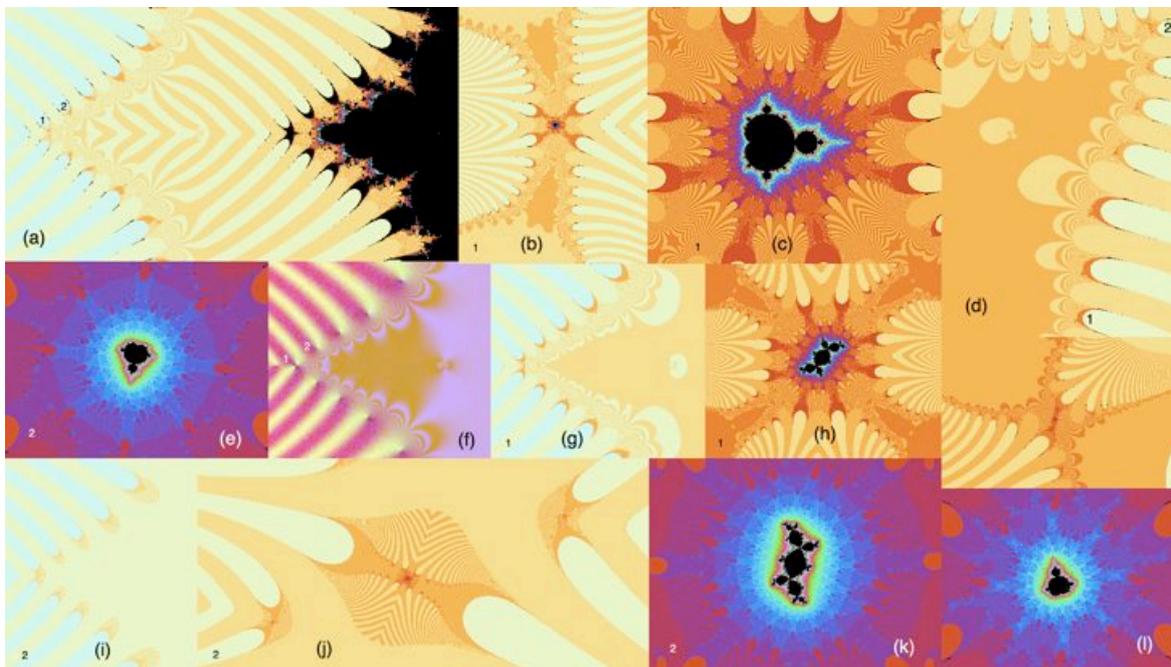

Fig 18: *z*-19 shows a further alpine displacement. The central valley is now displaced from a second valley – the 'principal valley' containing the principal point (1 f) and fixed values far to the left, (a). Its Mandelbrot set (b,c) is now at the horizontal fusion between successive conjoined frond ridges (b). The Julia set of a period 3 bulb of the principal Mandelbrot (g,h) shows homlogous structure. The fixed value at 2 in (a) points to a fractal recursion of sub-valleys at 1 in (d) rather than to the locus 2 in (d) where there is a Mandelbrot satellite in the same relative position as in the principal valley (e) with period-3 Julia (i,j,k). This continues with fractal replicates in successive 'unreal' fronds (l).

When we move on to *z*-19, the displacements have become even more acute. The transfer function now places the principal point and fixed values far into the negative, forming a shadow valley the

'principal valley' separate from the central valley. It is here we find the principal Mandelbrot set now nestled horizontally between two successive fused frond ridges, rather than vertically in a frond pair as previously. As before, his pattern is fractally replicated in the valley in each frond.

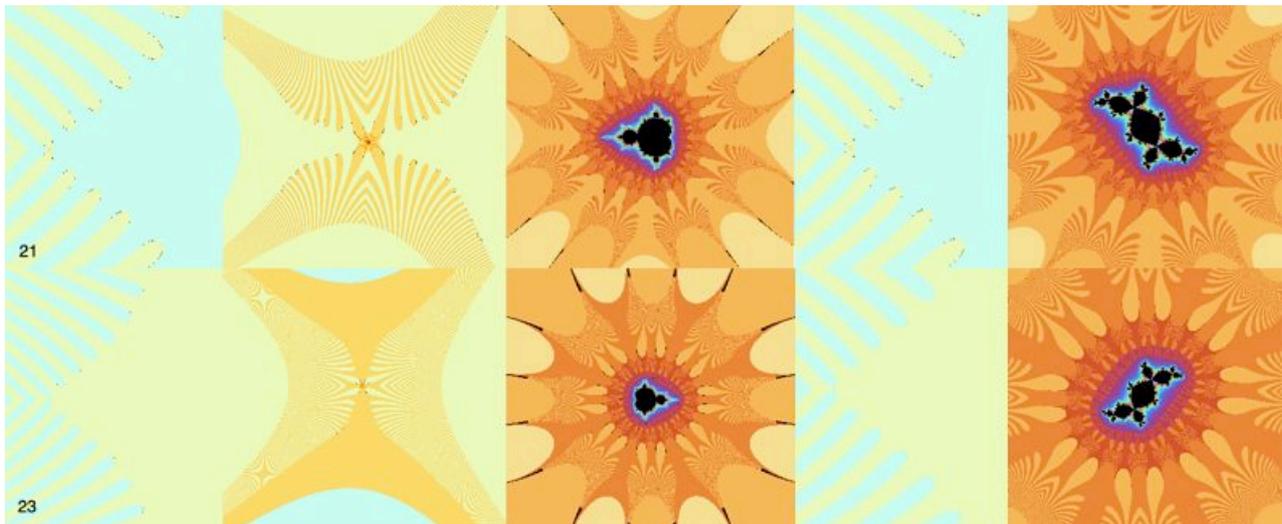
Fig 19: *z*-21 and *z*-23. The pattern of exponentiating maxima and minima corresponding to the series of fronds now continues with the maxima and minima following the structures of z-17 and z-19 displaced by ever huger positive and negative real translations.

The alternating pattern between *z*-17 and *z*-19 becomes a continuing sequence, evident in *z*-21 and *z*-23, where the central valley has now become entirely lost from view, enabling us to predict the dynamics of all subsequent criticals.

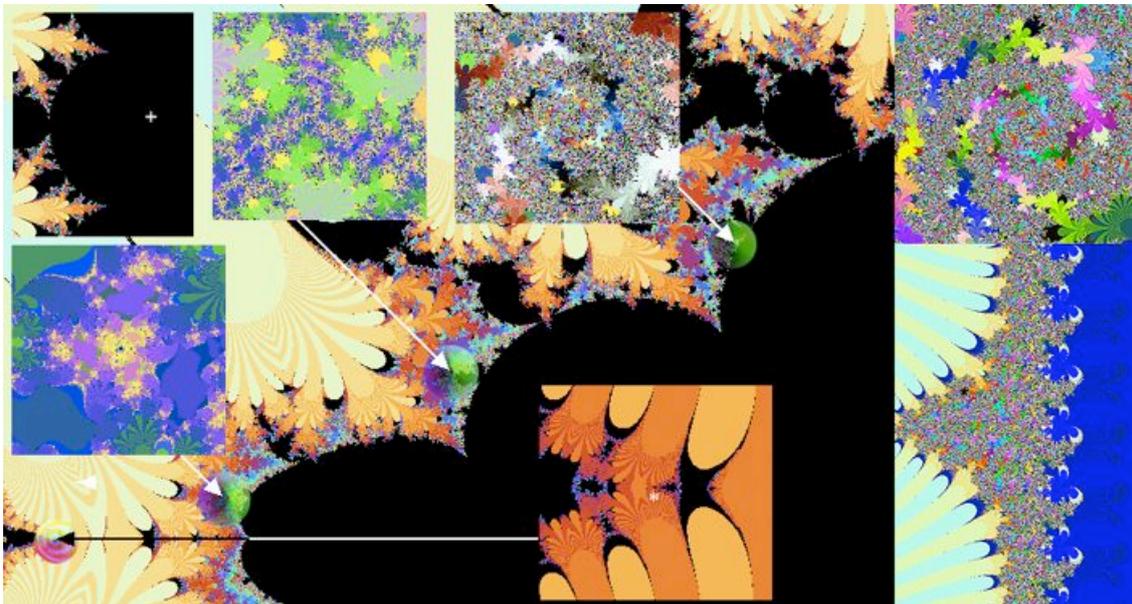
Fig 20: The unreal critical at z23 has its principal point well into the black ocean (+ top left), so we do not see distinctive features in its neighbourhood. The fixed values in the central valley lie outside the black ocean and correspond to two different types of M-points, the top two appearing as dendrite hubs and the other two are recursive fractal centres. For example the fixed value of the lower image (*) points to a valley at the base of a fractal valley ad infinitum. The top three all have fractal symmetries consistent with period 3. The Julia set of the top centre one (right) shows that this point is also an organizing centre of Julia dynamics in the neighbourhood of the fixed value (top right).

## (3) Shang-ri-La – The Unreal Criticals

We now turn our attention to the unreal criticals interspersed between the notorious non-trival zeros on the critical line x = ½, a little to the right of the zeros, with values from *x* ~ 0.78 - 2.4.

The locations of the critical points are generally to the right of the critical line and since their critical values are small their principal values lie close to the critical points in the Mandelbrot ocean. However some of them that are close enough to the chaotic landscape create local bays with quadratic Mandelbrot shorelines similar to the dynamics of $z$-15 in the central bay.

The first unreal critical z23 has a real part of $x = 2.4$ and shows little evidence of polynomial dynamics in the bay. All its fixed values in the central valley lie in chaotic territory, either at apparent dendrite hubs or loci of an endlessly recursive fractal process. The derivative function confirms these should all be repelling and thus constitute Misiurewicz points. All of those off the real axis appear to have period 3 symmetry. The Julia sets generated by these fixed values display a centre at the same fixed value with homologous dynamics.

The dynamics in the Julia set of fig 20 demonstrates that critical points far away from the central valley, can influence the dynamics there around their fixed values. The Julia parameter is simply the M-point corresponding to the fixed value on the boundary of the central valley, not the unreal critical point $z$23, yet the distant unreal critical is leaving its mark on the Julia set defined by a $c$ value in the central valley. This shows the 'Arizona effect' – the humming you hear on the telegraph wires out in the silent desert is a superposition of vibrations potentially coming all the way from California. In a similar way, the complex boundary of the central valley for each critical point is a combined 'whispering' of all the critical points, both real and unreal, which is why it is complex and sometimes highly amorphous. We shall see later that there are often Mandelbrot satellites in the neighbourhood of repelling fixed values, but for additive unreal criticals, these may all be submerged in the Mandelbrot ocean.

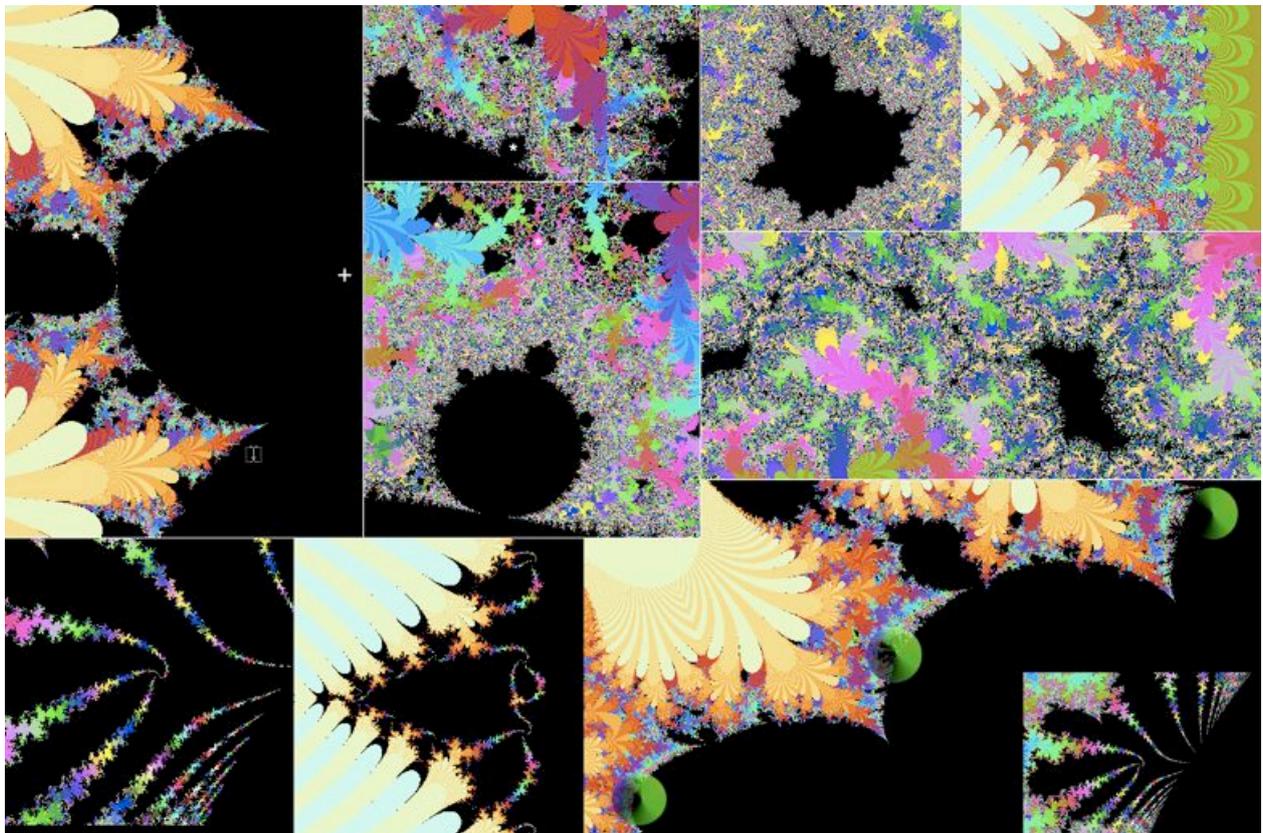

Fig 21: Dynamics of z31. (Upper sequence) shows the local basin of the critical point with quadratic bulbs, and a series of exploded views from the starred points to a Mandelbrot satellite whose Julia set has a (low resolution) period 3 kernel web. (Lower sequence) the central valley has two fixed values lying within the black ocean and only one (centre) on the boundary, pointing at a triple vertex of three clefts at the branching structure inset. The same structure on $z$-7 is a fractal version of the central bay of the same kind as the 'ant' of fig 10 and the branched pattern is also visible as distortions of the $z$-15 quadratic bulbs in the top of $z$95 in fig 11.

With *z*31, the second unreal critical, with a real value of 1.29, we begin to see richer polynomial dynamics. The bay bounding the region of the principal point now has a series of quadratic bulbs and these have dendrites supporting well-defined Mandelbrot satellites, which also give rise to period 3 Julia kernels from their period 3 bulbs as shown in fig 21. In this case two of the fixed values in the central valley lie in the ocean and only the centre one tends to a boundary M-point, this time at the triple vertex of three frond tips, (Mandelbrot cusps), with a Julia set having homologous dynamical centres.

As a third example, we have *z*95, which has a low real value of 0.78 and lies in a small focused bay with prominent quadratic bulbs, having dendrites supporting chains Mandelbrot satellites, whose period 3 bulbs generate confirmatory period 3 Julia kernels, establishing classic polynomial dynamics associated with an isolated critical point.

Again, this has only one of its fixed values in the central basin on shore, where it forms a fractal centre, again of period 3 nature. This suggests that much of the complex amorphous structure in the shoreline of the central bay is a product of the interaction of a large number of the unreal criticals with similar fixed value locations to those of *z*-15 acting together in superposition.

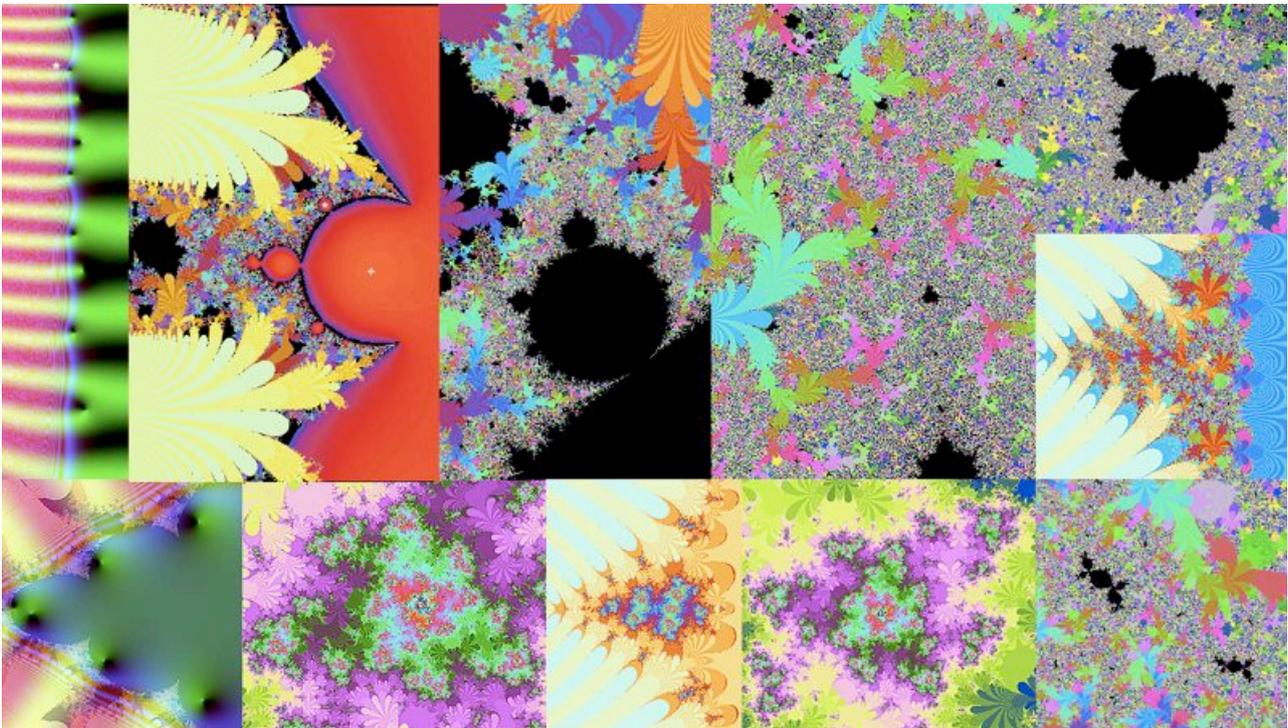

Fig 22: Dynamics of z95 (* top left). In the upper sequence is shown the location of z95=0.78+95.29i with critical value 0.43+0.078i with real part lower than that of *z*-15. The low real part of the critical point's coordinates causes it to be nestled closely towards the shore of the ocean, giving rise to a well-formed quadratic basin highlighted to show the iterations around the unreal critical (+). A series of exploded views from the bulb (*) leads to a well-formed satellite Mandelbrot, whose Julia set has well-defined web (lower right) of period 3 kernels. The lower sequence shows only the most left-hand fixed value lies on shore and gives rise to a fractal centre, again of period 3 symmetry, whose Julia set again has a homologous dynamic around the location of the fixed value.

**B: The Multiplicative Universe**

We now enter the universe of the multiplicative Mandelbrot parameter planes $\zeta_c(z) = c\zeta(z)$ a very different cosmos, where many of the patterns we have discovered to date will be manifest in new and different ways.

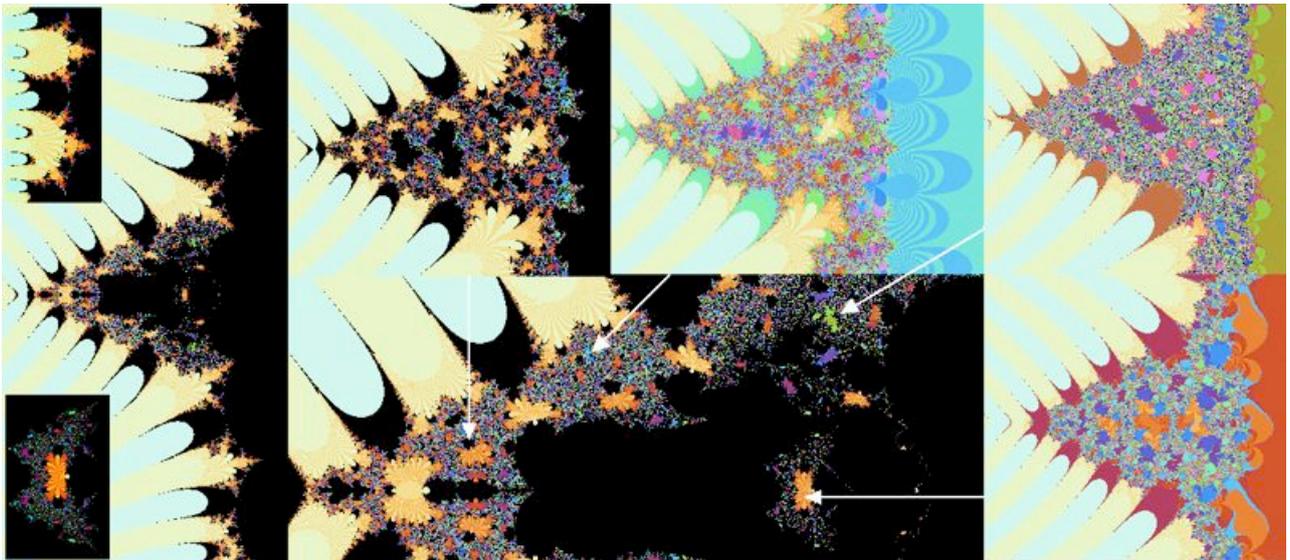

Fig 23: The global outlines of the multiplicative Mandelbrot set for the asymptotic plateau show fundamental similarities to the additive case with some notable differences. The large fronds are now phased with the critical points, rather than the zeta zeros which can now be found in the region of the blue bays. The central valley now has a thick chaotic exponential crown. As with the additive case (fig 5) points on this parameter plane classify asymptotic dynamics in the right half plane, leading to explosions during motion on curves in parameter space.

## (1) The Far Horizon

The first case we examine, is that of the asymptotic plateau, where, as with the additive case, the parameter plane forms an atlas of the dynamics on the asymptotic plateau as shown in fig 23. As previously, paths in the parameter plane $M_\infty$ can result in explosions in motion videos of Julia sets.

The additive and multiplicative parameter planes have superficial similarity, but these hide significant differences. The fronds are much larger and neatly phased around the critical points with the regular angular variations of dzeta, while the zeta zeros (and indeed all the real zeros from -18 on) are located in the zero point of the first order fractal valleys, since
$1000 \to z_n \zeta(1000) = z_n.1 = z_n \to z_n \zeta(z_n) = z_n.0 = 0$ and the fractal valleys are immediate pre-images

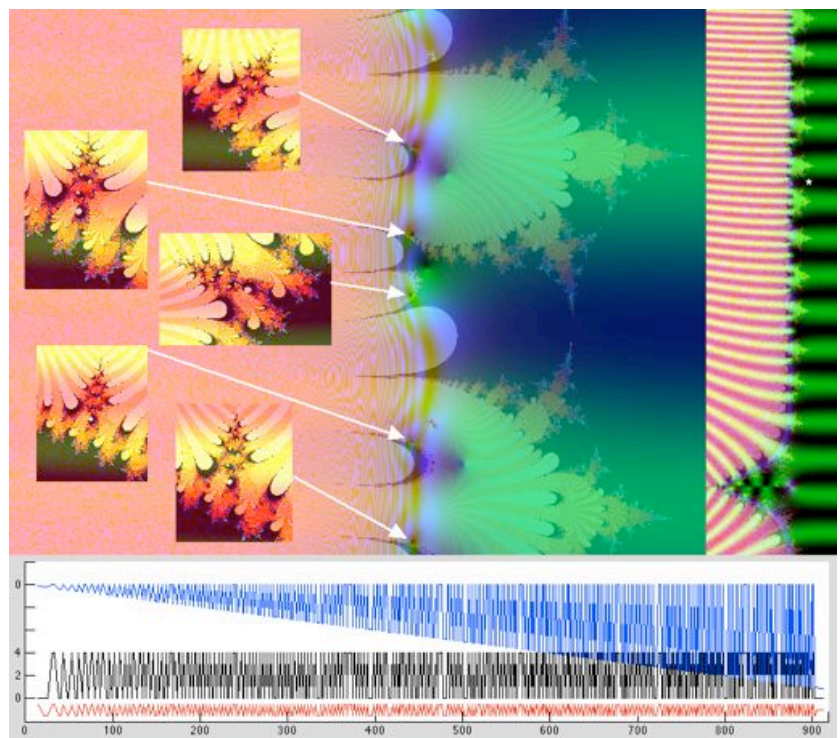

of the central valley surrounding zero, although the zeros and their neighbourhoods subsequently diverge, because of their large $c$ values. Hence each fractal valley is a dynamic map of the neighbourhood of each zeta zero.

Fig 24: The large fronds of $z1000$ are phased with the critical points of zeta (large dimples) and the regular angular variations of dzeta (right) while the zeros (small dimples) are located in fractal sub-bays on the shanks of the fronds. $z1000$ is ideally placed to model zeta itself because its critical value is 1 so there is no rescaling. Lower the first 600 non-trivial zeros showing intermittent zeros have periodic dynamics, iterating to themselves (blue) with period 4 (black) determined by the sign of the real part of the third iterate (red), while others diverge to $-\infty$.

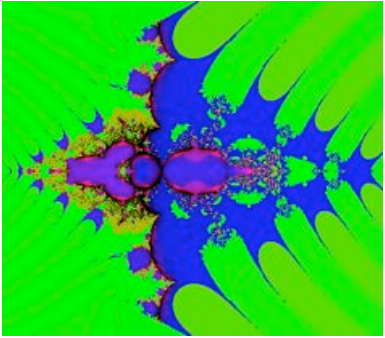

Fig 24b: $M_0$ gives an atlas of the dynamics of 0 for varying *c* values. (Fixed points and periods blue through red and yellow green. Escaping points green.)

The zeros iterate $0 \to z_n \zeta(0) = -0.5 z_n$, in the exponentiating left half-plane. Around half (286/600) iterate far into the right half-plane and enter strongly attracting period 4 cycles, since they are drawn back close to $z_n$ as with 1000, as in fig 24, e.g.
$z_{613} = 0.5 + 613.59i \to 0 \to -0.25 - 306.8i \to 353.98 + 11665i \to z_{613}$
Since the iteration passes through 0, the pseudo-Mandelbrot set $M_0$ from the non-critical 0 forms an atlas of the eventual dynamics of all the zeros for every value of *c*, $0 \to c\zeta(0) = -0.5c \to \cdots$ as in fig 24b, giving the zero dynamics for each multiplicative Julia set.

**(2) Tierra Firma – The Rift Valley**

If we examine the changes in the region around the central valley as we go through the real criticals, the results are surprising and unexpected. The entire picture of the right half plane being asymptotically fixed in the black ocean is lost and we find the exponentiating regions are now approaching from both positive and negative reals, with a rift valley in the centre with two opposing kinds of dynamics – on the left the central bay with submerged dynamics amid fractal replicates of the singularity island and on the right fractal structures of islands interspersed with Mandelbrot sets.

In fig 25 the rift valleys are shown for *z*-2 through to *z*-17. From there the rescaling takes the dynamic onto such small scales, the central valley ceases to exist. In fig 26 we examine the fractal structures in the opposing side of the rift valley for *z*-2 through to *z*-11, where the satellites form part of a lattice of islands whose symmetries are uniquely determined for each critical point.

Fig 25: The Rift Valley: As we move through the smaller real critical points, the central arena flexes wildly in tectonic tsunamis, with 'evolving' exponentiating regions (lighter yellow and ochre) overlaying and interacting with a 'static' resonance with the asymptotic limit's dynamics (darker ochre shades), folding around the complex plane, so that they also appear in the positive real half plane, giving rise to fractal structures for positive real values. Once again *z*-15 forms a boundary between the miniscule criticals and the vast criticals beyond.

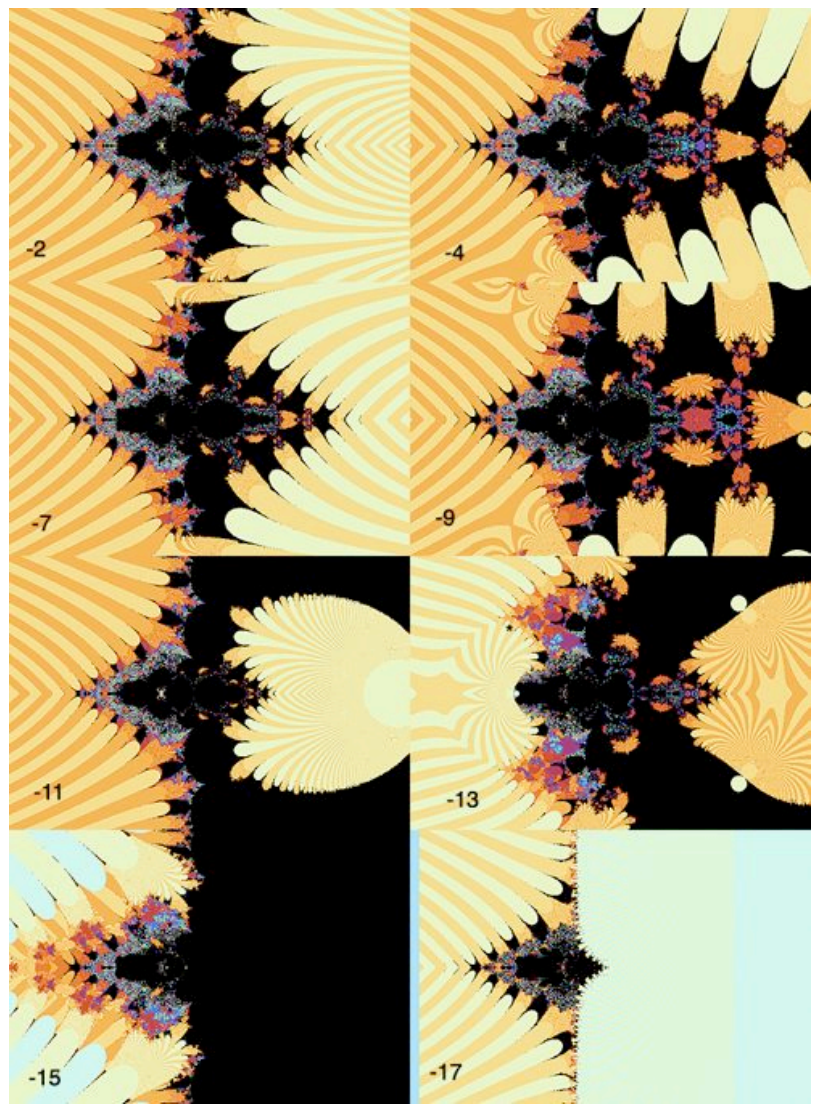

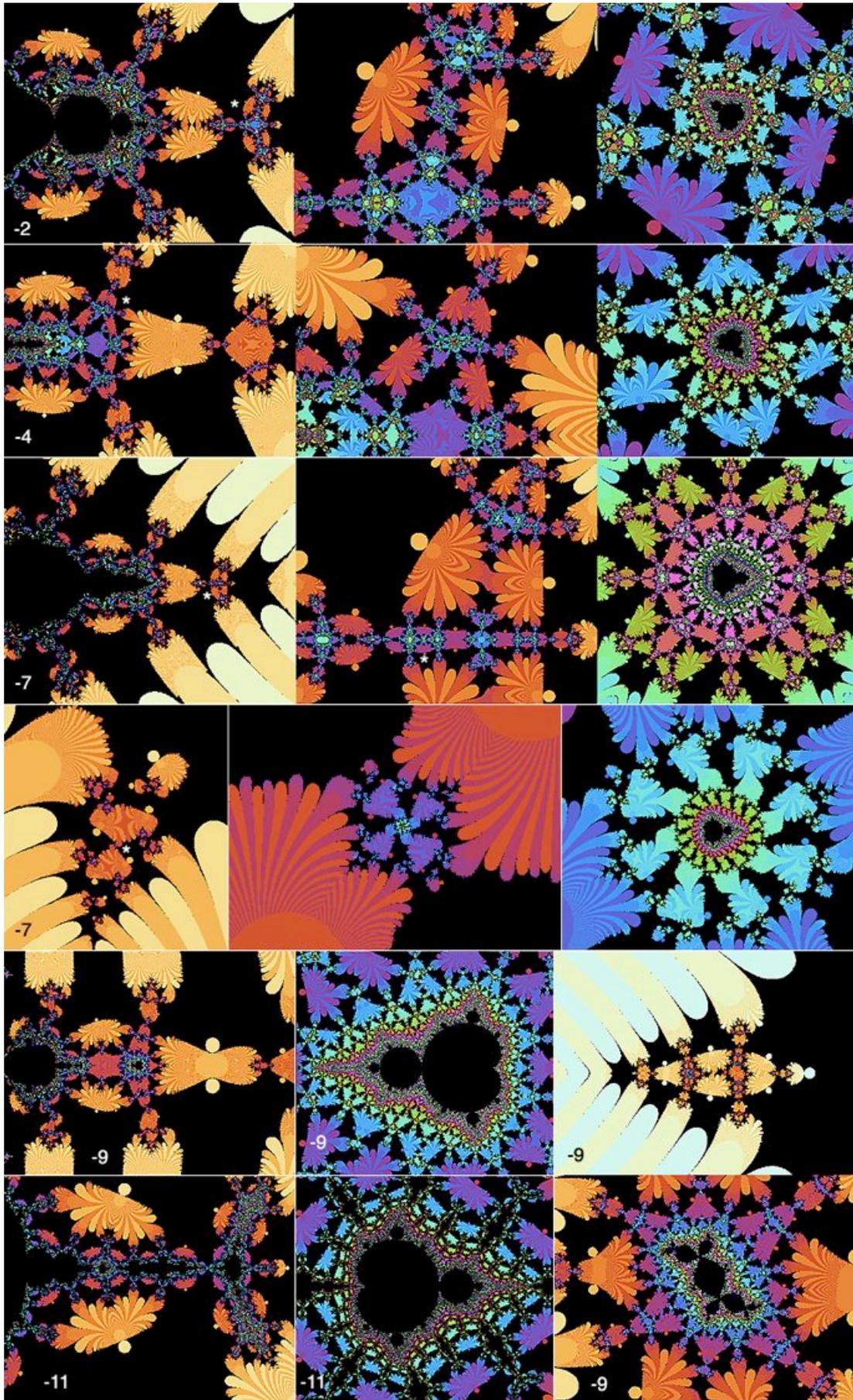

Fig 26: The positive arena of the rift valley displays a series of generic fractal lattices, supporting a network of Mandelbrot satellites (*) from *z*-2 to *z*-11 each with their own unique structure. Lower right the Julia set for *z*-9, which can be compared with that of *z*-2 in fig 30. Fractal versions of these structures also occur in fractal valleys in fronds above and below the rift valley in z-7 (4th row) and z-4.

# (3) Fractal Cosmology in the Galactic Abyss

We now turn to the multiplicative transfer functions to look for the principal points and fixed values of the miniscule criticals. The transfer functions take the form $\zeta(cv_r) - v_r = 0$ and the principal point is defined by $c_p : c_r = c_p \zeta(c_r) = c_r / \zeta(c_r) = c_r / v_r$. In both cases the process is being rescaled by being divided by the critical value, which for the miniscules varies down to $4 \times 10^{-3}$.

The principal valley of the transfer function is now huge – up to 20000 units across. Looking at $z$-2 to $z$-13 we find a consistent picture in which the principal valley alternates between large positive and large negative real values, ad the critical value changes sign, with a large principal Mandelbrot set and repelling fractal structures at the other fixed values, which support Mandelbrot satellites with higher periods within the local fractal structure, as illustrated in fig 27, where the rift valley has now become vanishingly small down the centre line of the figure.

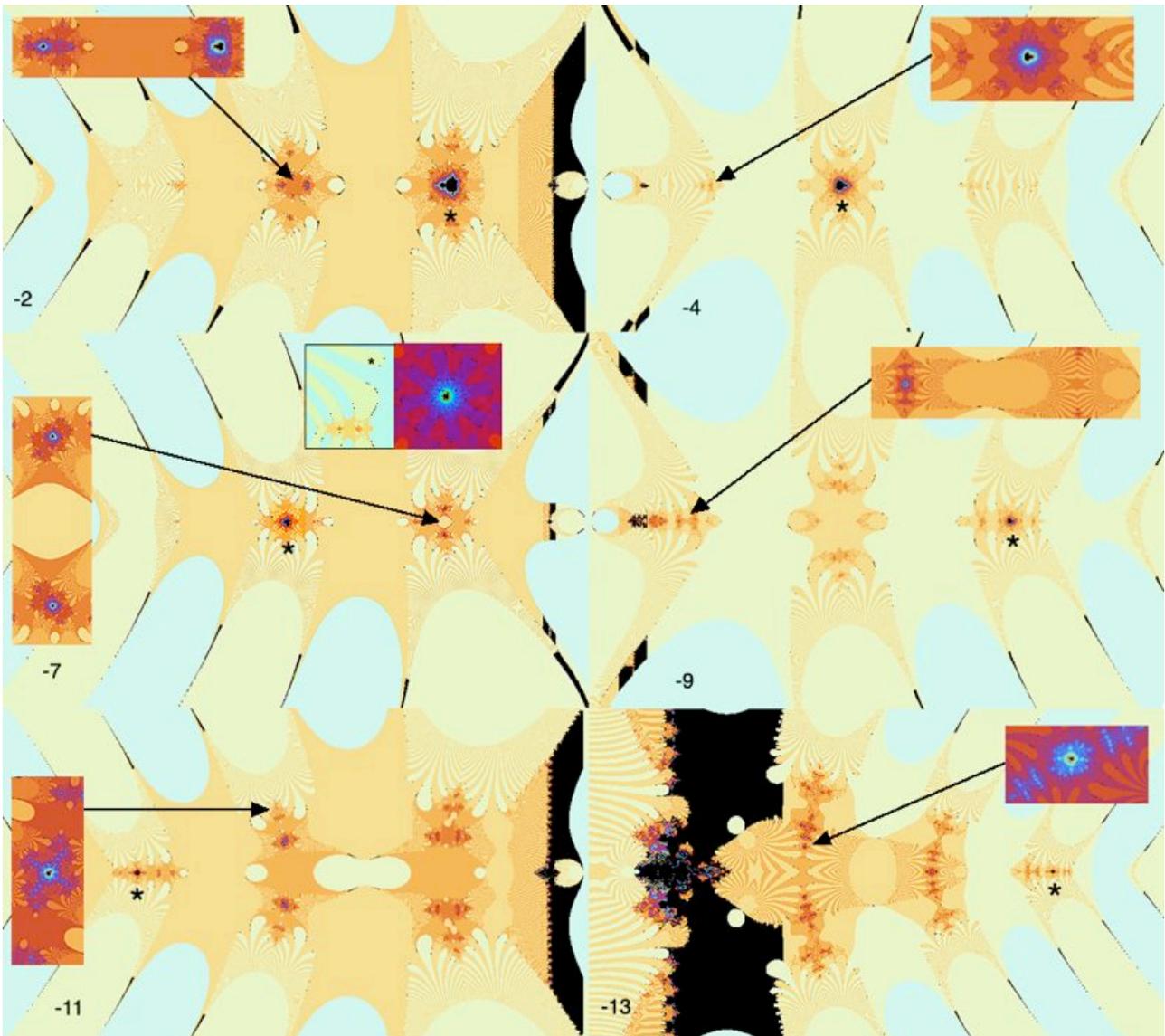

Fig 27: For the 'miniscule criticals', the principal valley defined by the transfer function is absolutely huge, reaching sizes of 2000, because the transfer function scales zeta by dividing it by the critical value. The central arenas can be seen tiny in the distance down the centre of the figure. Alternating critical maxima and minima flip principal valleys between positive and negative half-planes as their critical values change sign from positive to negative. The principal points each support a Mandelbrot set. As with the additive world, the locations follow a graded evolution (*) and like the vast criticals of the additive world, fixed values correspond to repelling M-points, which support fractal regions which contain Mandelbrot satellites. Inset in $z$-7 is a satellite from the first unreal valley above the principal valley.

Again, there is a progression across the fronds, with a maximum and a minimum in a progression down each frond accompanied by flipping between positive and negative half-planes, because the local maxima are positive and the minima are negative, causing negative scaling. In fig 28, we confirm by incipient period colouring that the satellites have higher periods, consistent with an attracting period, although the fixed value is repelling.

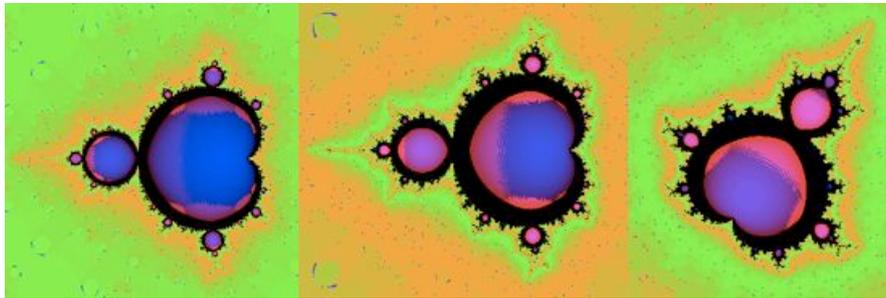

Fig 28: The principal Mandelbrot and two satellites of *z*-4 with colour-coded periods of 1, 2 and 3. Investigation of the periodicity of the Mandelbrot satellites confirms the satellites have higher periods, which enables them to be attracting even though dzeta may have absolute value greater than unity at their location.

**(4) Nanoworld's Myriad Valleys**

We can also explore the progression of the critical points in miniature, by examining fractal replicates of the principal valley in the regions to the left of the central valley, as illustrated in fig 29. The progression occurs in a similar pattern to the additive world as there is no flipping of signs.

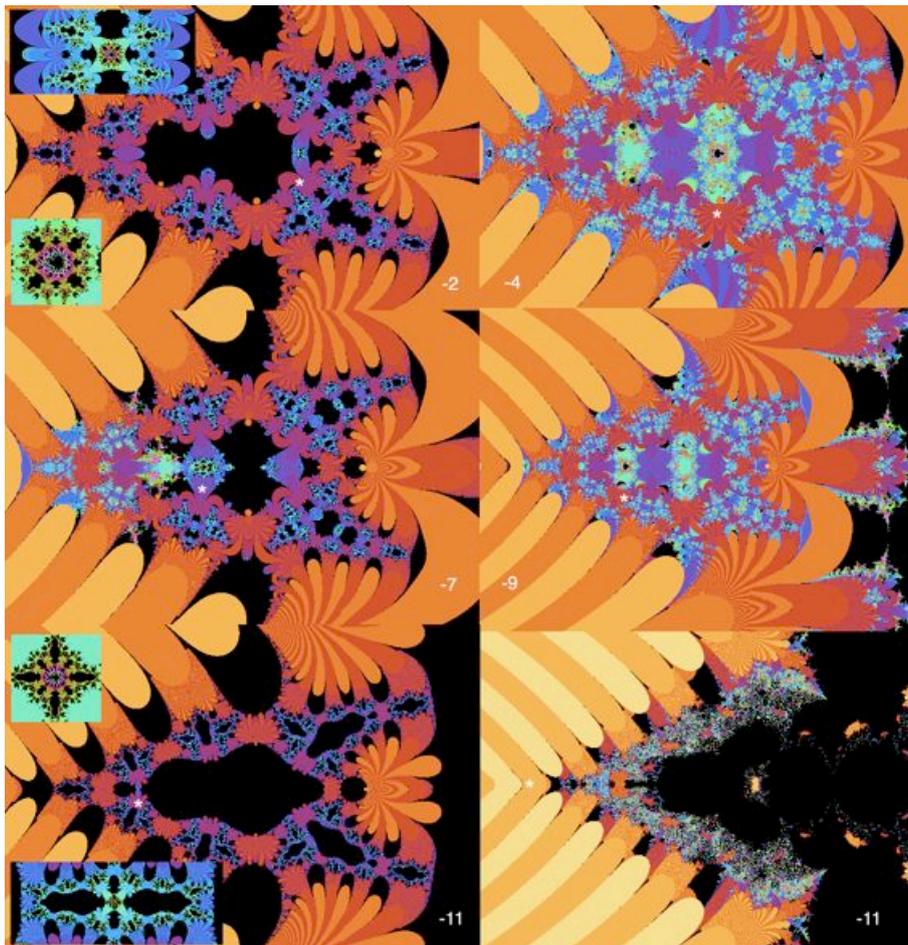

Fig 29: We can also examine the evolution of the 'miniscule criticals' in microscopic fractal replicates of the principal valley (*) lower right. The Mandelbrot satellite of the principal point (*) again evolves through the frond locations.

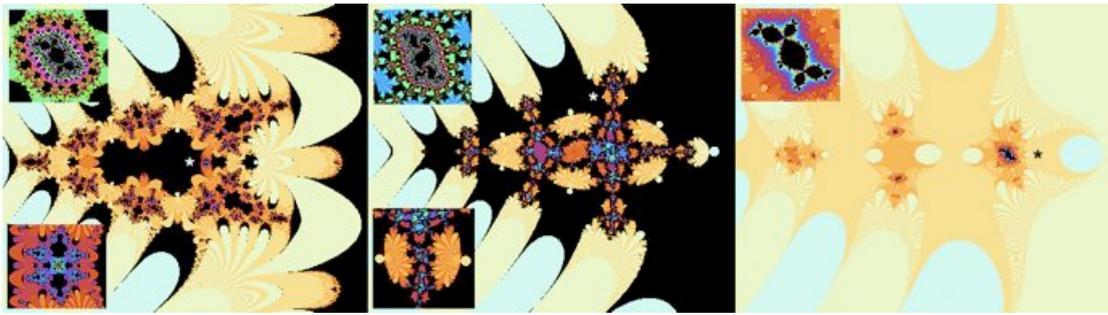

Fig 30: Three typical Julia sets of period 3 bulbs of Mandelbrot sets from the three regions - the microscopic cleft (left), the positive arena (centre) and the principal valley (right).

## (5) Jewels in the Crown

The critical point *z*-15 again forms a transition between the dynamics of the miniscules and the vast criticals. The positive half-plane has again become part of the main Mandelbrot set, but because the critical value is 0.52, the transfer function is scaled so the central valley region is scaled by ~2.

Fig 31: Again *z*-15 forms a transition between the miniscule and vast criticals. Because its critical value is 0.52, the principal valley size is doubled, so the principal point and its Mandelbrot lie at the extreme rear of the central valley (1). Again the fixed values are repelling and one can find satellite Mandelbrots located in fractal replicates of the central valley (2) and in successive fractal replicates of valleys in the fronds (3,4,5). There is also evidence for other repelling M-points of *z*-15 in the crown (far right). The case of *z*-13 remains somewhat enigmatic with no obvious Mandelbrot satellites in the positive arena and only in the two clefts surrounding the central valley (* in fig 24).

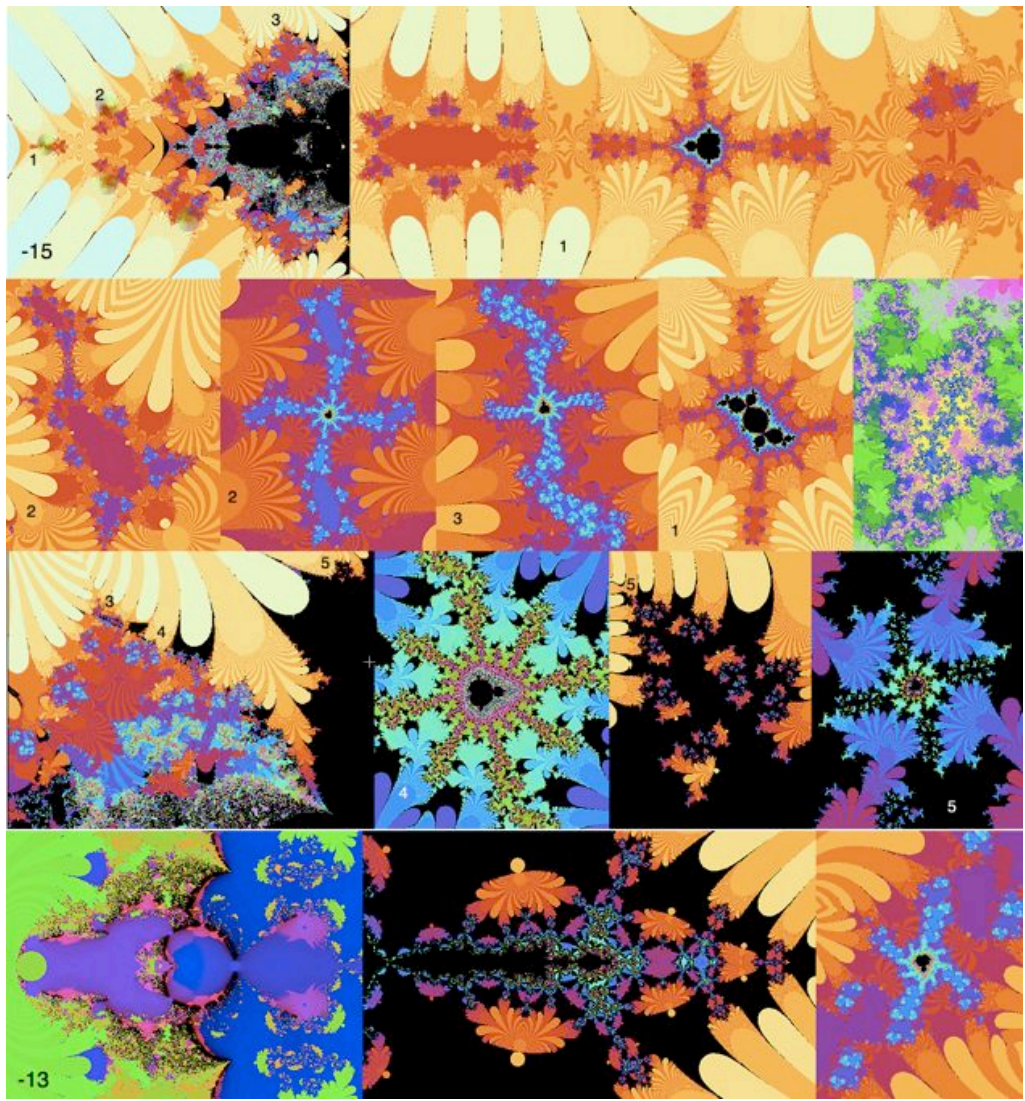

This results in the principal point and its Mandelbrot set being far back in the base of the bay and the fixed values being high up in the fronds of the crown. As with the miniscules, the fixed values are repelling, but satellite Mandelbrot sets are fractally repeated in the valley of each frond and in each generation of sub-valley as shown in fig 31.

*z*-13 also forms a transitional case, enigmatic for its lack of satellites in the positive arena of the rift valley and also lacking any fractal structure in the base of the central valley, due to the obstructing exponential node.

**(6) The Heavenly Heights**

When we move on to the vast criticals, we find a situation similar to that of the additive case although geometrically inverted in scale. Because the critical values are now exponentially large, the principal valley is successively downscaled into a minute region near zero. At *z*-17 we find the central valley opposed by a copy of itself very similar to the additive case of *z*-17, again enclosing the Mandelbrot set of the principal point vertically in a pair of fronds as in fig 32, fractally replicated in sub-bays as shown right.

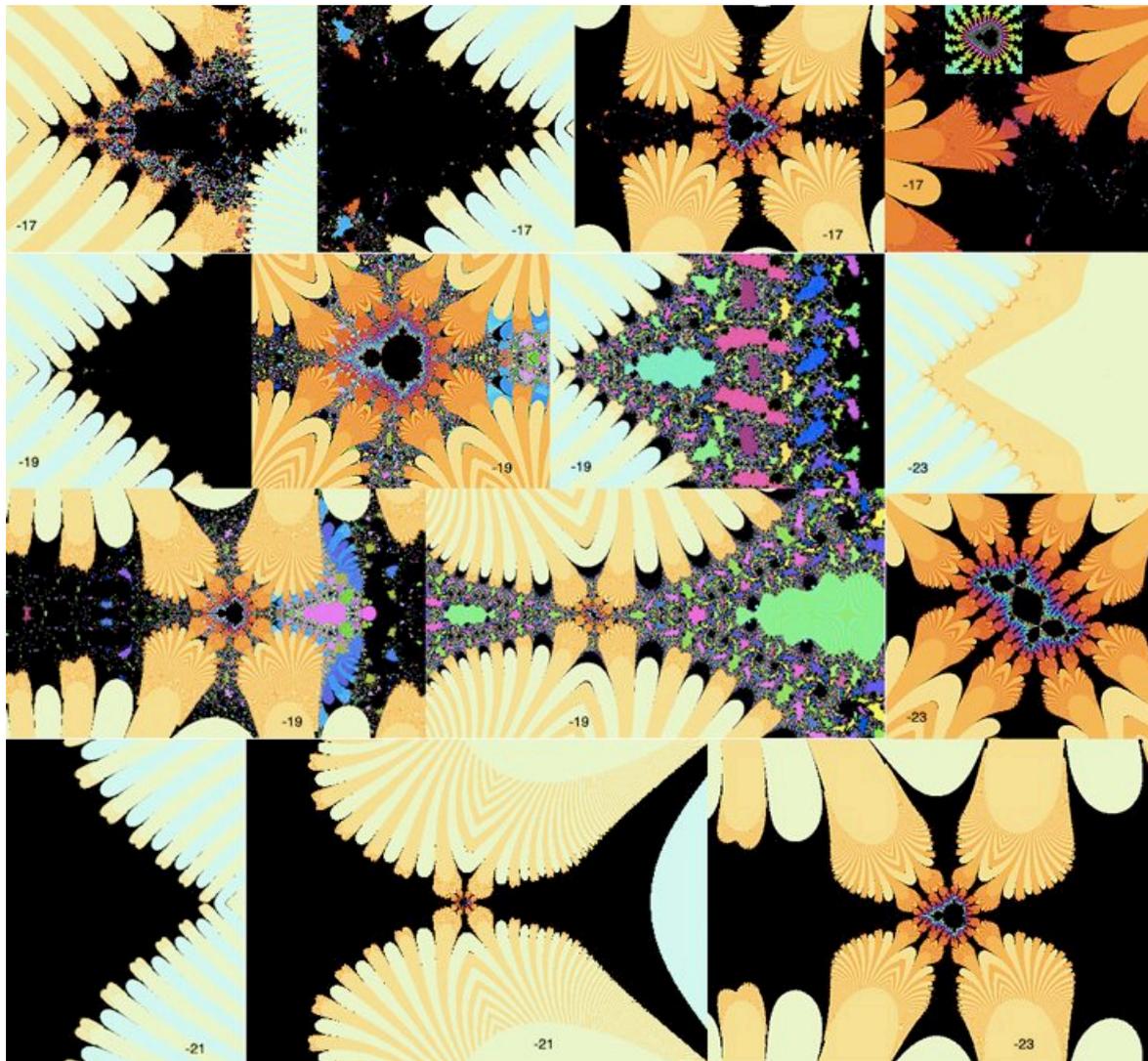

Fig 32: The 'vast' criticals again follow an evolution, but slightly different from the additive scenario. The central valleys have now become microscopic because the critical values are huge and alternate between positive and negative half-planes according to the sign of the critical value, as with the 'miniscule' criticals. In all cases the Mandelbrot sets are pinched vertically between fronds, and as before the process is repeated in the fractal valleys of each frond.

As with the miniscules, successive vast criticals have principal points flipping between positive and negative half-planes, because their critical values change sign. Fig 32 shows the evolution for the first four vast criticals, including Julia verification for *z*-19 and *z*-23. Rather than alternating between vertical and horizontal frond placement, all principal Mandelbrot sets are held vertically, but their locations flip between positive and negative half-planes, but on miniscule scales, so that their Julia sets are generated covertly from *c* values near zero.

**(7) Pandora's Dystopia – The Dissonance of the Unreal Criticals**

To cap the bag, we have the principal Mandelbrot sets of a series of unreal criticals up to $z95$. Rather than being in the Mandelbrot ocean, the principal sets can end up in a variety of locations, oceanic, far into the chaotic exponential region, or nestled in the tips of fronds. Some have irregular amorphous forms due to dynamic interference. Their locations can be far into the positive or negative due to the twisting of the transfer function by the complex number critical values. Significantly, as with the real criticals, in figs 27, 31 and 32, one can find satellite Mandelbrot sets in the neighbourhood of (repelling) fixed values in neighbouring fronds to the principal point.

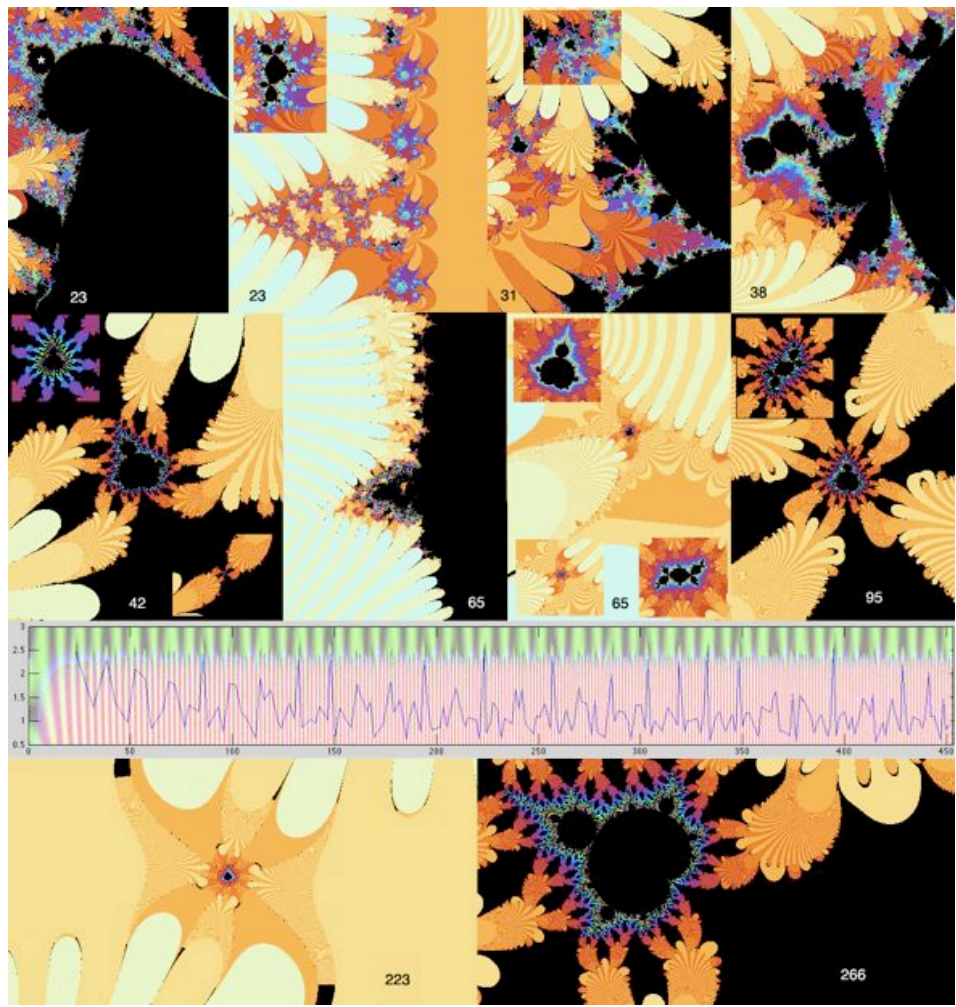

Fig 33: We finally come to the multiplicative unreal criticals where chaos bites back. Their transfer functions can be tilted, (see $z65$) through division by the critical value so the principal points can occur with large positive and negative real values. The first ($z23$) lies in the Mandelbrot ocean and has quadratic bulbs at the shoreline. Some have perfectly formed Mandelbrot sets around their principal points, which may lie deep within the exponential chaotic region or be delicately supported by the tips of fronds. However in several cases the principal sets are highly irregular. Well-formed period-3 Julia sets can be gained both from the bulbs of quadratic basins ($z23$) and Mandelbrot sets ($z65$) but the period-3 Julia kernel is readily detectable only in the vicinity of the critical point where it is largest, rather than in a disseminated web, so can easily be missed. With $z95$ = 0.78+95.29i, the primary Mandelbrot is at the principal point 40.7+241.71i and the Julia set has this as *c* value, but we need to examine it at the critical point to find the period 3 Julia kernel. It takes more terms of the zeta sum as we explore greater imaginary values, at last 500 for $z95$ to present without serious distortions of the Mandelbrot set, which becomes computationally time-consuming. Lower strip, a graph of the real and imaginary positions of the zeta criticals superimposed on the zeros of dzeta in the range up to gamma underflow to 0, ~450i. To take matters to the limit, we include principal sets for $z223$ = 2.500042+223.408567i and $z266$ = 0.791525+266.014624i, to similar scale, whose principal points lie at 19.696+ 232.45i and -217.38+475.88i respectively, the latter reaching to imaginary values on the border of underflow of the gamma function, and requiring up to 1500 zeta terms on 'Full'-terms iteration. Finally, if we were right when we saw with real criticals, that fixed values, although often repelling, had Mandelbrot satellites in their neighbourhoods, we should be able to find them for unreal criticals as well. The inset in $z31$ with principal point 1.8190+44.8408i shows one instance, in a similar position in the lower frond at ~4.044+34.484i, neighbouring the adjacent fixed value, implying the same pattern pertains. In $z42$ insets is also shown a satellite at ~8.794+46.720i, near the adjacent fixed value, while the principal point of the main set is at 10.0451+51.7145i.

The Julia sets generated by *c* values the period 3 bulbs of the principal sets and satellites all correctly display period 3 kernels, however they have an enigmatic quality, unlike previous examples, where the period 3 kernel is disseminated through the Julia set. In this case the kernels are strongly centered on the critical point location (see *z23* in fig 33) rather than the central valley or principal point, so they can easily be missed. Accurate calculation for these critical values and principal points requires much higher numbers of zeta terms.

**Appendix: Fractal Geography of Eta, Xi and Dirichlet *L*-functions**

As noted in the introduction, zeta is defined through the Dirichlet Eta function:

$\eta(z) = \sum_{n=1}^{\infty}(-1)^{n+1}n^{-z}$ which shares the same trivial and non-trivial zeros but also has regular zeros where $1 - 2^{1-z} = 0, z = 1 \pm 2k\pi i$. Riemann also introduced a symmetric function Xi, which also has the non-trivial zeroes $\xi(z) = \Gamma\left(\frac{z}{2}+1\right)(z-1)\pi^{-\frac{z}{2}}\zeta(z)$. Because of their close relationship with zeta a summary of their fractal geographies is included for comparison.

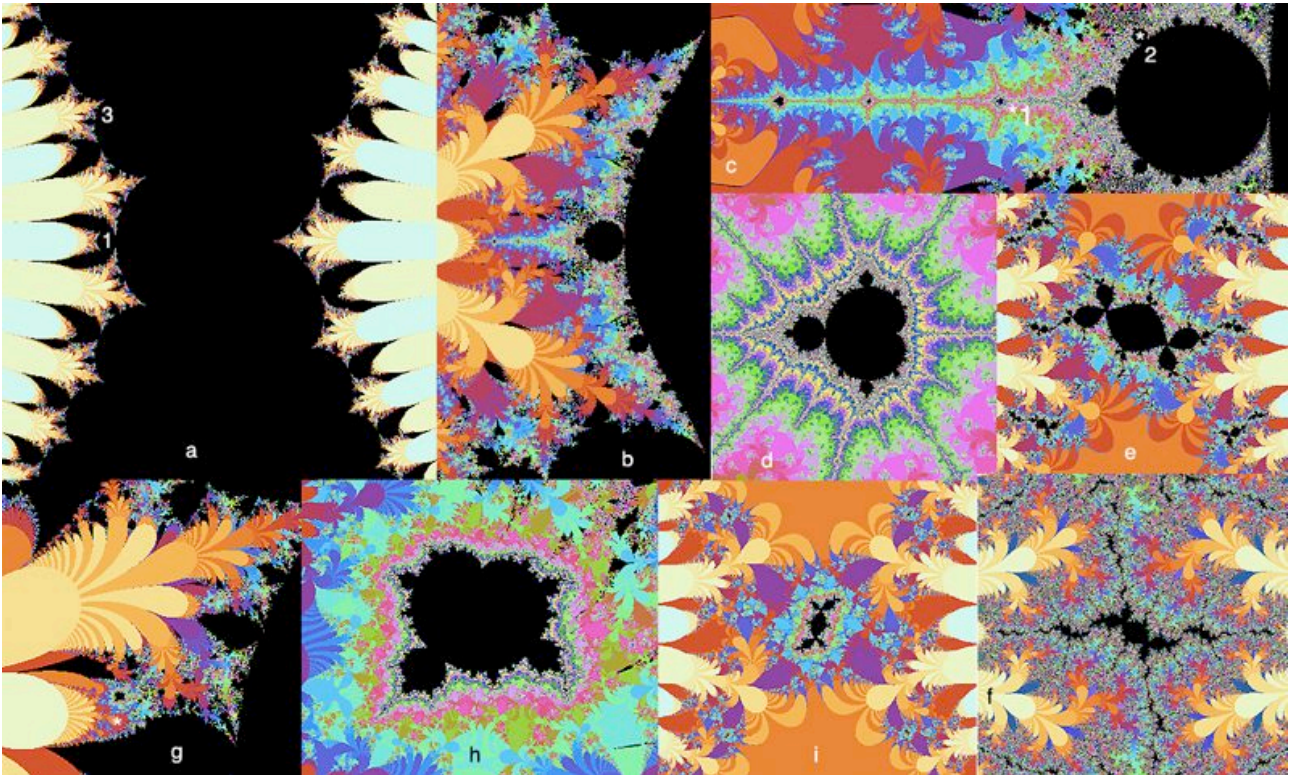

Fig a1: Additive Mandelbrot set of xi from *z*=½ (a) has a small centre-left frond (a1, b) with multiple quadratic bulbs (c) supporting satellites (d). Their period 3 bulbs generate period 3 Julia sets (e, f). Other fronds (a3, g) also have satellites (h) correctly generating Julia dynamics (i).

Xi has critical points running up the critical line, but the function is very close to 0 throughout, so only the central critical at *z*=½ leaves any mark on the Mandelbrot dynamics. This results in a very simple profile for Xi similar to that of our initial example $f(z) = ze^{-z}$, with local quadratic dynamics and global exponential dynamics in the additive case, and only exponential fronds in the multiplicative case, which is omitted.

Eta has a more intriguing geography, because it shares the zeta zeros, and also has an additional regular set of unreal zeros, resulting in a more widely scattered set of unreal criticals. Its real criticals form a shorter sequence, with much clearer manifestations of quadratic dynamics throughout. Its unreal criticals are even more varied, although similar qualitatively to zeta.

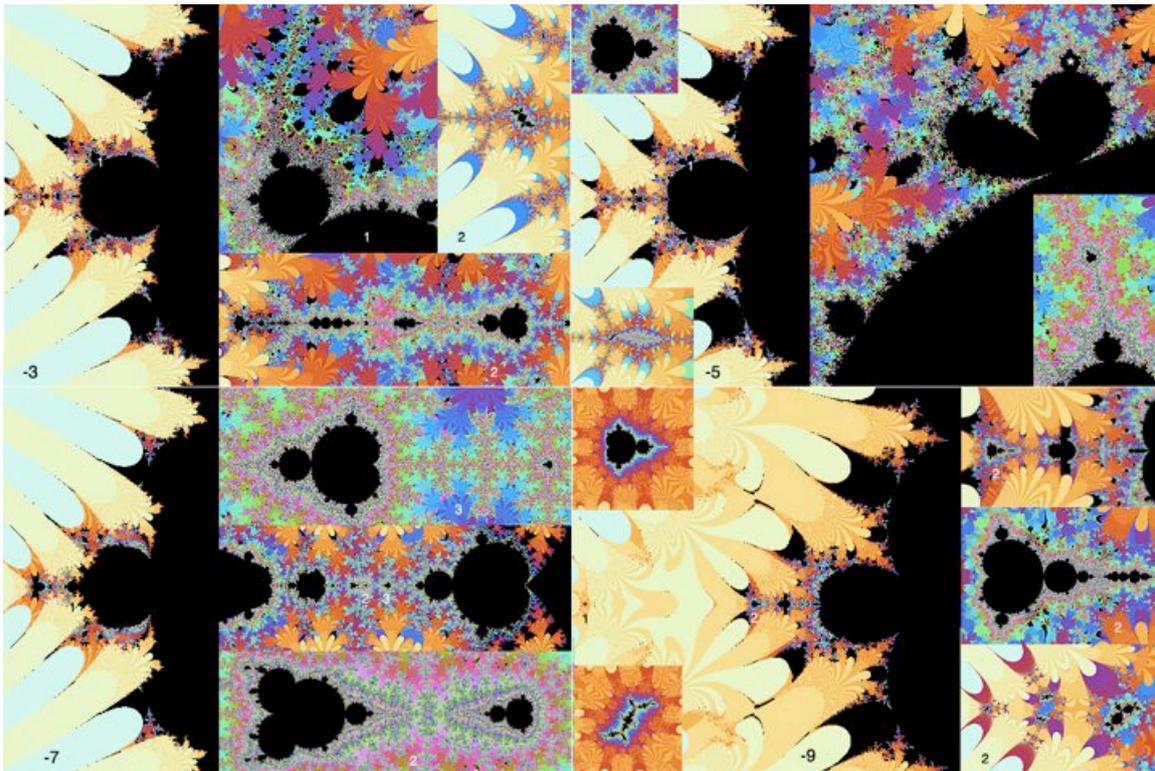

Fig a2: Additive Mandelbrot and Julia dynamics for the first four real criticals of eta. These follow a similar trend to those of zeta, with *e*-3 having bulbs in the outer basin, *e*-5 bulbs on the inner, *e*-7 having a Mandelbrot held between the inner two fronds and *e*-9 having a principal Mandelbrot far into the negative (1). However all four also have Mandelbrot satellites supporting Julia dynamics along the real axis and from the bulbs in *e*-3 and *e*-5.

The details of the fractal geometry of eta's dynamics is summarized in the figures and captions.

Fig a3: While the unreal critical that have larger positive real parts are too far into the Mandelbrot ocean to display quadratic dynamics on the shoreline, eta's unreal criticals whose real parts are less than 1 form additive quadratic basins between their local fronds, with satellites generating corresponding Julia dynamics. Upper sequence *e*72. Lower sequence *e*33 and *e*99.

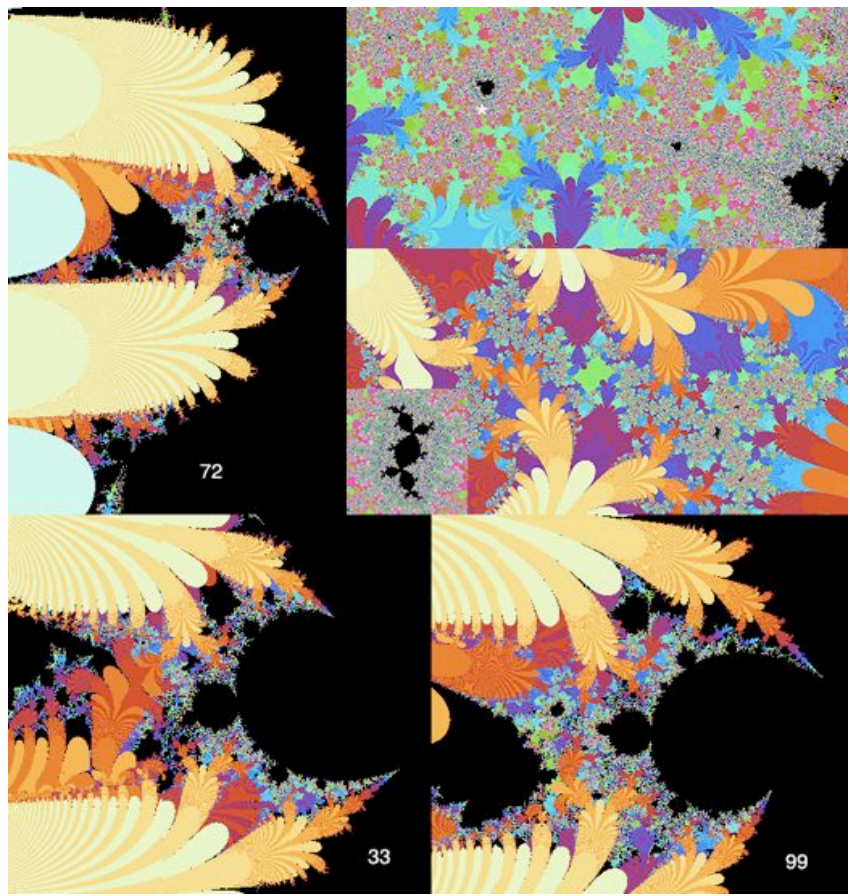

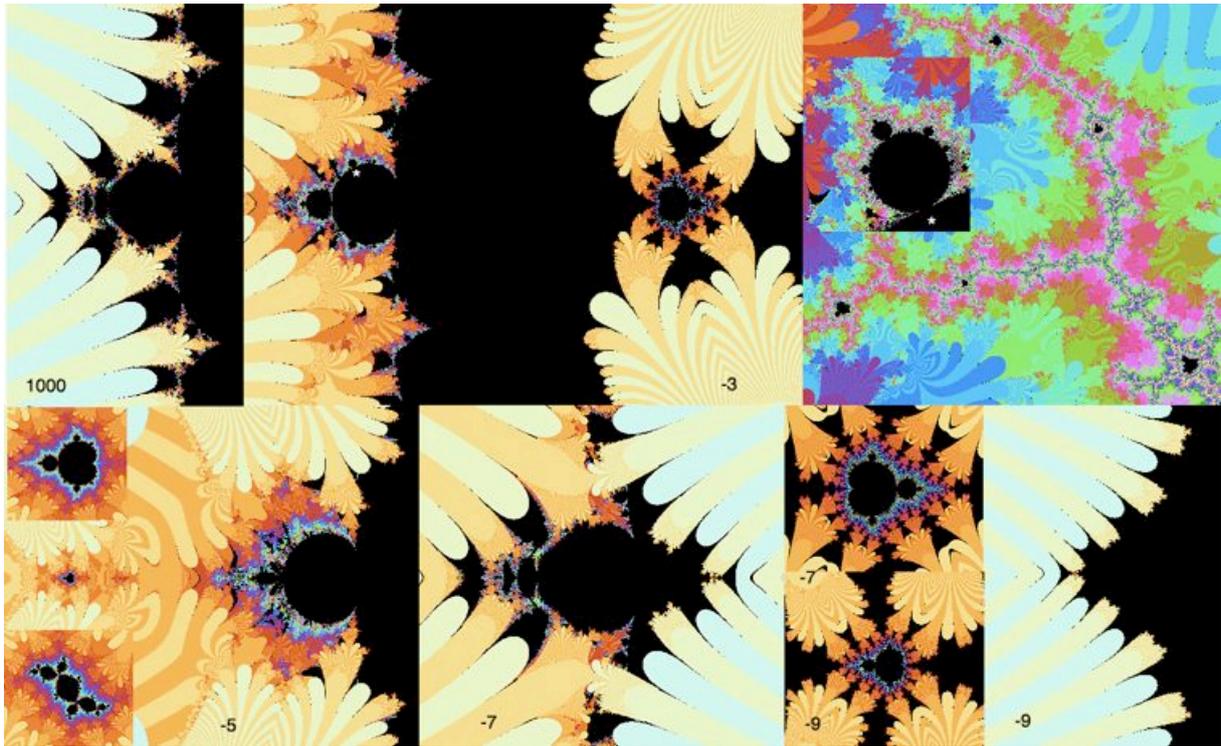

Fig a4: Eta's multiplicative Mandelbrot sets show a similar evolution to zeta with important differences. *e*-3 has both a principal Mandelbrot set in the positive half-plane and unlike the chaotic crown of fractal islands in zeta, a central bay with quadratic bulbs and satellites. *e*-5 has a principal set in the negative half-plane and *e*-7 and *e*-9 have principal sets on alternating sides. There is no sign of the series of lattice structures zeta had in the positive arena of the rift valley. The dynamics of *e*-7 in the central valley shows a strong correspondence with *e*1000.

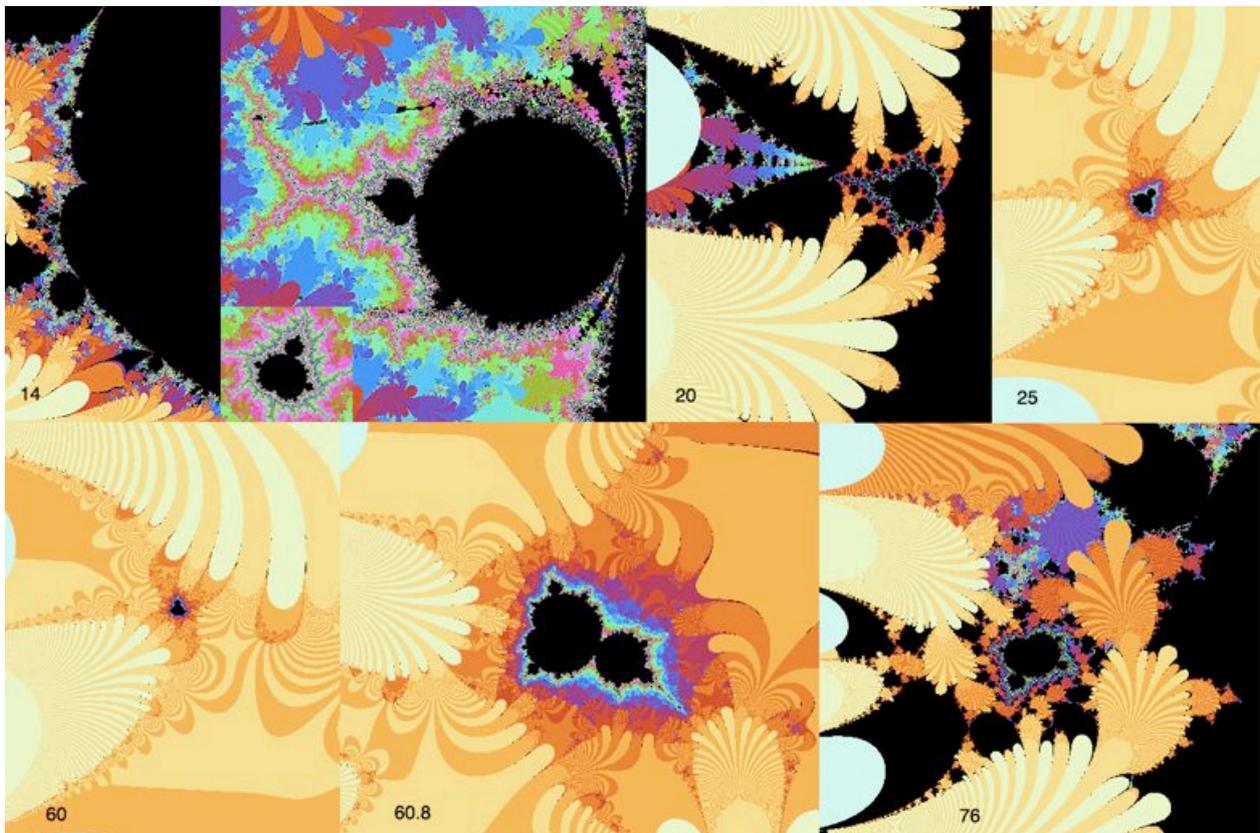

Fig a5: Eta's unreal criticals form a varied sequence of cases, from the first one (*e*14), with a large real part lying in the Mandelbrot ocean with quadratic bulbs supporting satellites, through the next (*e*20) touching tips of frond dendrites to the next (*e*25) interior to the chaotic domain. Some critical points (*e*60, *e*60.8) are scattered laterally (1.09 and 2.42), one of which has a much larger real value but their principal sets are still in the chaotic escaping region. In the most extreme case examined (*e*76), with a real part of 2.52, we are just back to a location enfolded in the frond tips.

The Dirichlet L-functions[6] $\sum_{n=1}^{\infty} \chi(n)n^{-z} = \prod_{p \text{ prime}} (1 - \chi(p)p^{-z})^{-1}$ where χ are a cyclic set of Dirichlet characters generated by a finite residue group, display several new properties of the fractal geography of zeta functions., some have a more complex twist, resulting in twin conjugate functions, which are asymmetric about the *x*-axis, as is *L*(5,2) in fig a6. While there are no degenerate zeros with multiple roots in zeta, as it can be expressed $\zeta(z) = \pi^{z/2} \frac{\prod_{\rho \text{ non-triv } 0} (1 - z/\rho)}{2(z-1)\Gamma(1 + z/2)}$, showing all the zeros are of first order, including those derived from gamma's singularities, the principal characters of numbers with *k*-distinct prime factors have a multiple zero at the origin with a neighbouring critical point having degenerate higher degree Mandelbrot sets, as noted in fig a6.

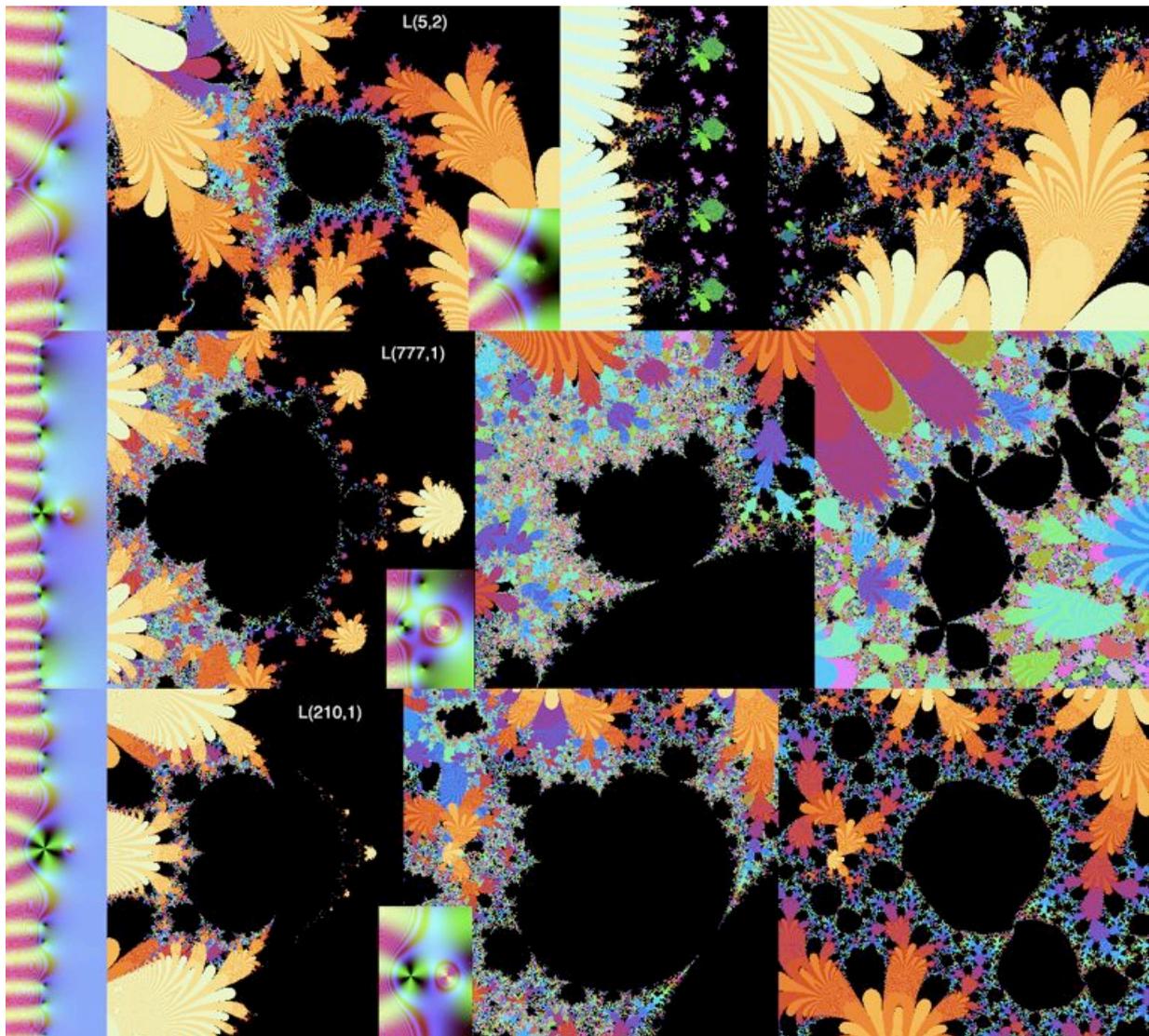

Fig a6: Function profiles, Mandelbrot and sample Julia sets from *L*(5,2), *L*(777,1) and *L*(210,1) show several manifestations of *L*-function fractal structure. The central *L*(5,2) Mandelbrot set is twisted asymmetrically because the coefficients exist in two conjugate forms. *L*(777,1) and *L*(210,1) show examples of degenerate cubic and quartic Mandelbrot sets due to the multiple degenerate critical points showing as zeros in the derivative (inset). These derive from the fact that each is a product of distinct primes - 777=3.7.39 and 210=2.3.5.7, resulting in factoring zeta by 3 and 4 product terms $(1-p^{-z})^{-1}$ respectively.

Continuing into the realm of abstract *L*-functions of elliptic and other curves (), we find new fractal properties emerging from the genus of the embedded curves, in which curves of higher genus give rise to degenerate higher degree Mandelbrot and Julia sets throughout their central valleys.

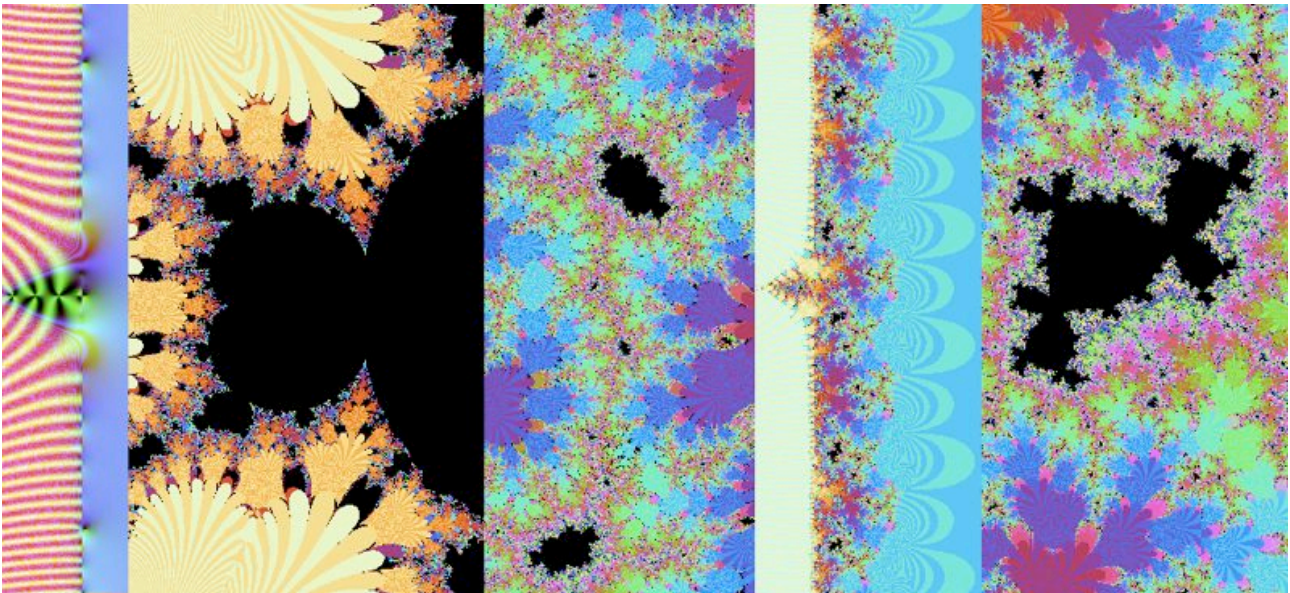

Fig a7: *L*-function of the genus-3 curve $y^2 + (x^3 + x^2 + x + 1)y = x^7 + 2x^6 + 2x^5 + x^4$ has an entire set of degenerate cubic zeros arising from its multiple repeated gamma factors. Its critical points correspondingly give rise to cubic Mandelbrot satellites with cubic Julia sets.

## Introduction for Non-mathematicians:

The aim of this paper is to present a tour of the zeta function which anyone with some passing understanding of maths can enjoy and appreciate, as an exciting intrepid journey on complex math space. The paper is based on an easy to use Mac application developed by the author using XCode, downloadable from: http://www.dhushara.com/DarkHeart/. You can thus all experience this journey yourselves and even make video records of your adventures, using the RZViewer application.

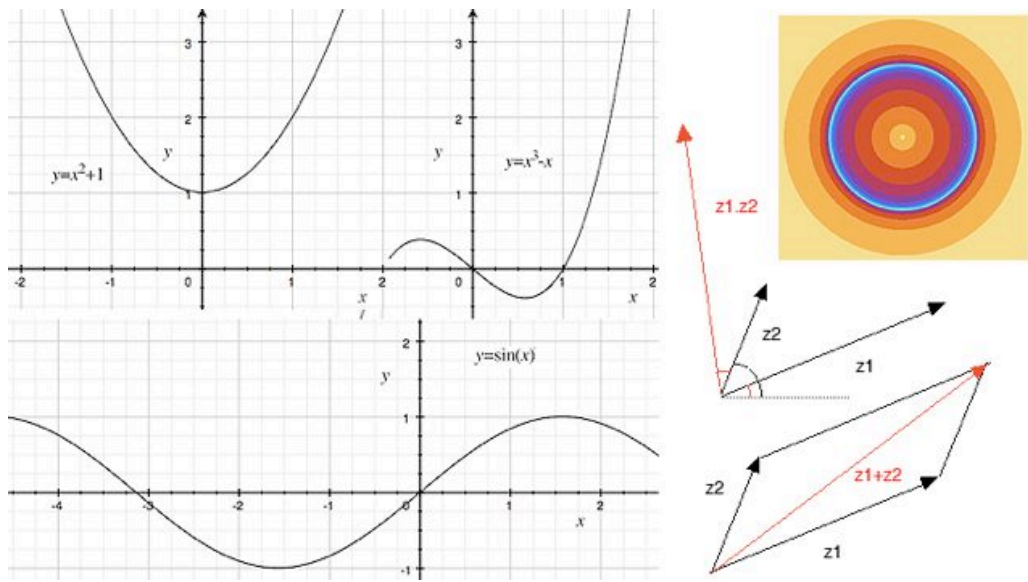

Fig 34: Critical points of real functions, complex addition and multiplication and the simplest Julia set.

**Complex Numbers, Fractals and Chaos**

Just about everyone has a good idea of the real number line and understands how they add and multiply, and how functions such as $y(x) = x^2$ can be depicted in graphs with humps and troughs. Functions with several terms in powers of *x*, such as $y(x) = x^3 - x$, are the polynomials. Virtually every other function we commonly use, including the so-called 'transcendental' trigonomeyric and exponential functions can be written out as an infinite polynomial in a power series: e.g.

$\sin(x) = x - \frac{x^3}{3!} + \frac{x^5}{5!} - \frac{x^7}{7!} + \cdots + \frac{x^n}{n!} + \cdots$, commonly written as $\sin(x) = \sum_{n=1}^{\infty} \frac{(-1)^n x^{2n+1}}{(2n+1)!}$. Sin(x) behaves precisely as a polynomial with infinitely many humps and troughs.

Complex numbers are more abstruse for many people, who find their use of imaginary numbers perplexing. We make numbers complex by trying to solve $x^2 + 1 = 0$, which has no zeros on the real line, by adding $i = \sqrt{-1}$ to the number system, forming the unit along a new y-axis in a 2-dimensional number plane. These numbers turn out not just to solve one equation, but provide existence of the zeros for all polynomial equations. Moreover they become essential in the study of wave processes in physics. Complex numbers add real and imaginary parts, just like real numbers, but multiplication has a new twist to it. The sizes $r$ of the numbers multiply, but their angles $\theta$ add: $z_1 \cdot z_2 = (r_1, \theta_1)(r_2, \theta_2) = (r_1 \cdot r_2, \theta_1 + \theta_2)$. This causes complex numbers and complex functions to twist as well as expand and contract, so $z^2$ is $r^2$ the size but at double the angle, as in fig 34.

Enter chaos and fractals. One of the most basic ways of representing dynamic change is to map a function into itself in a feedback loop that causes a kind of stroboscopic 'discrete' flow in a series of jumps of iteration of the feedback loop $x_{n+1} = f(x_n)$. When we do this, the discrete flow falls into two kinds of dynamics - order and chaos - each of which tends to have opposing regions in the space in which the process is taking place. In a regime of order, the discrete iteration steps down a basin of attraction converging towards a defined fixed equilibrium, or periodic oscillation, rather like a stroboscopic picture of the water flowing out of a bathtub. In the chaotic regime, the discrete flow is divergent, so that arbitrarily close trajectories diverge exponentially –commonly known as the butterfly effect – arbitrarily small perturbations leading to global instabilities.

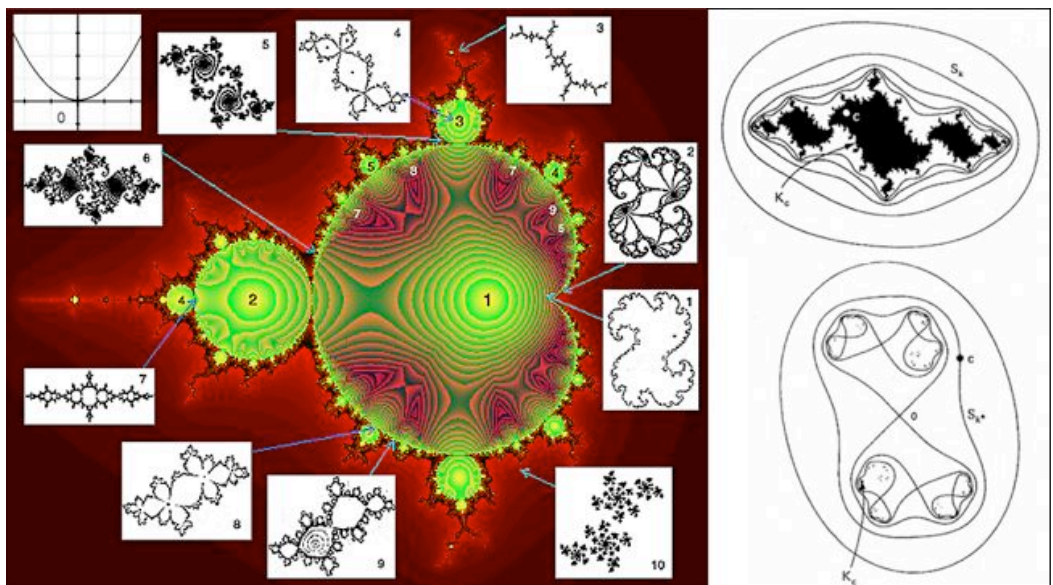

Fig 35: (Left) The Mandelbrot set of $f(z) = z^2 + c$ surrounded by its Julia sets. (Right) How the critical point of the function determines whether the Julia set is connected or not. The critical point $x = 0$ (see upper left) is the last point to escape. When the Julia set is connected (upper right) the 'filled-in' Julia set $K_c$ (black) contains the critical point and critical value $f(0) = c$. The inverse images of a large circle are all closed curves and the complement of $K_c$ is a disc. In the disconnected case the inverse images bifurcate to form a figure 8 and the circles repeatedly double. The centre of the first figure 8 is the critical point and the critical value c escapes to infinity. The Mandelbrot set is formed by iterating the critical point for each c value thus forming a 'connectivity atlas" of the Julia sets.

Discrete chaos happens in real and complex numbers alike, but in complex numbers the process becomes fully fledged, due to the twisting action of complex numbers and in the 2D plane of complex numbers, the dynamics becomes apparent in a way we can all experience.

**Complex Discrete Dynamics**
In 1918 Gaston Julia devised a method for iterating rational functions of complex numbers to produce a discrete dynamical system by feeding back a function again and again into itself. A rational function is a fraction of two polynomials so it's a little more general than polynomials like

$f(z) = z^2 + c$. The two regions of order and chaos became the complementary Fatou and Julia sets, and over time Julia sets proved to be fractals – geometric sets of points in which a process is endlessly replicated recursively in the manner of a snowflake.

The simplest Julia set is that of $f(z) = z^2$ simply the unit circle illustrated top right in fig 34. For numbers with absolute value greater than 1, repeated squaring makes them bigger and bigger and they race off to infinity. For smaller numbers squaring makes them smaller and they converge to zero. But for numbers on the unit circle squaring keeps doubling the angle without changing the size, so each point races round the circle chaotically at a different exponentiating rate. We can even see the butterfly effect of chaos because any two numbers have the difference of their angles doubled for each squaring step, so no matter how close they start out, after a finite number of steps they will get to be half a revolution apart and be on opposite sides of the circle.

Then Benoit Mandelbrot, who had invented the term fractal, began to study the set of points *c* in the plane, for which the Julia sets were topologically connected. This is a supremely subtle and difficult question, but it turned out to have an elegantly simple computational solution. For a quadratic, such as $f_c(z) = z^2 + c$, the Julia sets consist either of a totally disconnected Cantor set [1], with all points outside escaping to infinity, or a connected fractal, with points inside drawn into their own basins. The critical points, lying horizontally on the humps or troughs of the function, are always the last points to escape, so if the Julia set is connected, the critical point must be trapped inside. In this simplest case, the quadratic has just one trough at the origin, where the derivative, or slope, is zero. So if we iterate for each *c* from the critical point and test whether it does escape in a finite number of steps, we can colour the *c* value of the Mandelbrot set (traditionally black) if the critical point didn't manage to escape and thus portray it, colouring the complement by how many steps it did take to escape.

The Mandelbrot set and its complement is thus a fractal atlas of the dynamics of each and every Julia set - a fractal whose fractal regions are infinitely varied, leading to the claim of it being the most complex and beautiful mathematical object in existence. The heart shape also shows the effect of multiplicative twisting of complex numbers, with the cusp, and each of the bulbs, displaying different forms of fractional rotation highlighted by the number of dendrites emerging from each bulb. The 'M-points', Misiurewicz points - tips and intersection hubs of dendrites if they are repelling, and the roots of the Mandelbrot set and its satellites if they are attracting, form key points in our investigation, which are accessible because they are eventually fixed or eventually periodic in a finite number of steps, as opposed to the asymptotic dynamics of points in a basin of attraction, or chaotically wandering. Fractional rotation periodicities of the bulbs and dendrites also follow mediants $\frac{a}{b} + \frac{c}{d} = \frac{a+c}{b+d}$ so, between period *b* and period *d* bulbs, is a smaller period *b+d* bulb.

The Mandelbrot set, as complex dynamic atlas generalizes to polynomials, rational functions, and even transcendental functions, such as trigonometric and exponential functions, all of which have discrete dynamics, which can be described in terms of parameter planes forming an atlas of their Julia sets, with generic correspondence to the quadratic case, representing a simple hump, or trough, corresponding to an isolated zero, and occasional degenerate higher dimensional cases. In rational and transcendental functions, the Julia set kernels form an infinite web, and in the transcendental case, rather than chaos and the Julia set existing only on the boundary of the infinite and finite basins, the closure of the entire infinite basin is chaotic, but otherwise the classification of the forms of the kernels of the Julia sets appears to be universal.

---

[1] Georg Cantor described the simplest case of a fractal: a unit interval in which the open middle third is removed recursively from each remaining interval *ad infinitum*.

## The Zeta Function

As noted at the beginning of the paper, the Zeta function[7] is defined in two radically different ways:

$$\zeta(z) = \sum_{n=1}^{\infty} n^{-z} = \prod_{p \text{ prime}} \left(1 - p^{-z}\right)^{-1} \quad [1]$$

The relationship between the sum formula over integers *n* and the product over primes *p* was discovered by Leonhard Euler and is equivalent to prime sieving. This has led to a deep connection between zeta and the theory of prime numbers, and in particular to the, as yet unproven, and possibly unprovable, yet ostensibly true, Riemann hypothesis - that all the 'non-trivial' zeros of zeta (i.e. those not on the negative real axis) lie on the critical line $x = \frac{1}{2}$.

While both sides of [1] are convergent only for real(*s*) > 1, we can see that we can easily extend convergence to real(*s*) > 0 using Dirichlet's Eta[8], by defining zeta as an alternating series:

$$\eta(z) = \left(1 - 2^{1-z}\right)\zeta(z) = \sum_{n=1}^{\infty}(-1)^{n+1}n^{-z}, \text{ so } \zeta(z) = \left(1 - 2^{1-z}\right)^{-1}\sum_{n=1}^{\infty}(-1)^{n+1}n^{-z}$$

This takes us all the way across the critical strip, including the non-trivial zeros.

Riemann[9] analytically extended the zeta function to the entire complex plane, except the simple infinity at *z* = 1, by considering the integral definition of the gamma function $\Gamma(z) = \int_0^{\infty} t^{z-1}e^{-t}dt$, leading to the formula $\zeta(z) = 2^z \pi^{-1+z} \sin\left(\frac{\pi z}{2}\right)\Gamma(1-z)\zeta(1-z), z < 0$, where $\Gamma(z)$, the generalization of the factorial function, is also extended to the negative half-plane by the analytic continuation $\Gamma(1-z) = \frac{\pi}{\sin(\pi z)\Gamma(z)}$ (although this is not needed to define zeta). A full derivation of zeta is included in the companion paper [10]. The trivial zeros arise from those of the sine on the negative real axis, and the non-trivial, from $\zeta(z)$ itself in the critical strip $0 \leq x \leq 1$. We now have the full picture in fig 1, with exponential growth outside the central valley in the centre.

Zeta is famous, not just for its supreme complexity, but its enigmatic non-trivial zeros, which all appear to lie symmetrically on the line $x = \frac{1}{2}$. Proving this has remained elusive, because their locations arise from an infinite number of 'holographic' superpositions of the powers of *n*, which vary exponentially in size with *x* but have sinusoidal wave functions of *y* varying logarithmically with *n*. Despite the apparent symmetry in the *x* direction, the dynamics is determined by both the real and imaginary parts acting together, and the imaginary parts form an irregular series of values, consisting of wave transformations of the irregular distribution of the primes.